\documentclass[11pt]{article}
\usepackage{amssymb}
\usepackage{amsmath,amsthm,amsfonts,amscd,bezier}
\usepackage{amssymb}
\usepackage[usenames]{color}
\usepackage[T1]{fontenc}       
\usepackage[all]{xy}
\usepackage[normalem]{ulem}

\newtheorem{lemma}{Lemma}[section]

\newtheorem{theorem}[lemma]{Theorem}
\newtheorem{remark}[lemma]{Remark}
\newtheorem{proposition}[lemma]{Proposition}
\newtheorem{corollary}[lemma]{Corollary}
\newtheorem{definition}[lemma]{Definition}
\newtheorem{example}[lemma]{Example}
\newcommand{\Dem}{\noindent{\sc Proof:\ \ }}
\newcommand{\cqd}{{\hfill $\rule{2mm}{2mm}$}\vspace{1cm}}

\def\C{\mathbb{C}}
\def\A{{\mathcal A}}
\def\O{{\mathcal O}}
\def\N{\mathbb{N}}

\date{ \ }

\begin{document}

\title{On the analytic equivalence of branches in $(n+1)$-space}

\author{Marcelo Escudeiro Hernandes, Welinton Anderson Rocha\\ and Maria Elenice Rodrigues Hernandes  
\thanks{The first author is partially supported by CNPq grant 304594/2024-5, the second author was supported by CAPES, and the authors are partially supported by CNPq grant 407454/2023-3.}
\\ 
{\small Departamento de Matem\'atica, Universidade Estadual de Maring\'a}\\
{\small Av. Colombo, 5790, 87020-900, Maring\'a - PR, Brazil \footnote{Email address: mehernandes@uem.br, welintonrocha1998@gmail.com and merhernandes@uem.br (corresponding author)
} }}

\maketitle

\begin{abstract}
In this paper, we consider Newton–Puiseux parametrizations of irreducible curves in $\mathbb{C}^{n+1}$, $n \geq 1$, within a fixed semigroup under the action of Mather's group $\mathcal{A}$. We establish criteria for eliminating parameters while preserving the Newton–Puiseux form, extending known results for plane curves.
\end{abstract}

\vspace{0.3cm}
\begin{center}
 2020 Mathematics Subject Classification: 14H20 (primary), 32S05, 32S10 (secondary)

\vspace{0.2cm}
Keywords: parametrized curves, semigroup, analytical equivalence.
\end{center}

\vspace{0.5cm}

\section{Introduction}

The aim of this paper is to investigate local properties of irreducible analytic curves (branches) in $\mathbb{C}^{n+1}$, or equivalently, their parametrizations $\phi: (\mathbb{C},0) \to (\mathbb{C}^{n+1},0)$ with respect to $\mathcal{A}$-equivalence. A central question within this framework is how normal forms can be obtained for this class of objects.

This work is inspired by the extensive literature on the study of plane curves. The topological classification of plane curves (as immersed germs) has received significant contributions from several mathematicians, such as Zariski, Brauner, and Burau. For more details, see \cite{BK}, as well as \cite{top-um} for the irreducible case and \cite{top-varios} for plane curves with several branches. In \cite{fernandes}, Fernandes shows that the topological equivalence coincides with the bi-Lipschitz classification for plane curves. Regarding the analytic classification with a fixed topological type, as proposed by Zariski, the problem was solved by Hefez and Hernandes in \cite{HH2} for plane branches and, more recently, by Hernandes and Rodrigues Hernandes in \cite{RodriguesHernandes} for reduced plane curves. It is well known that the value semigroup, or equivalent data, constitutes a powerful tool for describing normal forms within an equivalence class. This is due to the fact that the value semigroup is a complete invariant for topological equivalence, that is, two plane curves are topologically equivalent if and only if their semigroups coincide. Since analytically equivalent plane branches are, in particular, topologically equivalent, the strategy adopted in \cite{HH2} consisted of fixing a topological class, namely, curves with the same semigroup, and distinguishing them by means of finer analytic invariants.

Several works have been dedicated to obtaining normal forms for curves in $\mathbb{C}^{n+1}$ with $n>1$. In \cite{Ebey}, Ebey presents a classification of algebroid curves under some particular conditions, which includes parametrizations with $\delta$ invariant less than $7$ (see Section 3 for definition of $\delta$). In \cite{GibsonHobbs}, Gibson and Hobbs classify the $\mathcal{A}$-simple parametrized curves in $\mathbb{C}^3$. Furthermore, Kolgushkin and Sadykov, in \cite{KS}, classify stably simple reduced curve singularities in a complex space of any dimension. Zhitomirskii \cite{Zh} classifies fully simple singularities of plane and space curves in the real case with any number of components, while Stevens \cite{ST} classifies simple parametrizations of complex curve singularities of arbitrary embedding dimension.

Although the value semigroup can be defined for branches in $\mathbb{C}^{n+1}$, it is not a topological invariant in this higher-dimensional setting; in fact, not even the multiplicity of a branch is preserved topologically (see \cite{lipman}). However, since the analytic equivalence of irreducible curves translates into the $\mathcal{A}$-equivalence of their parametrizations, it can be verified that the value semigroup remains an $\mathcal{A}$-invariant for branches (see \cite{GibsonHobbs}). The focus of this work is the study of parametrizations of irreducible curves in $\mathbb{C}^{n+1}$ within a fixed semigroup, under the action of the group $\mathcal{A}$, in order to obtain normal forms.

Considering a primitive parametrized curve in $\mathbb{C}^{n+1}$, we recall that any branch admits an $\mathcal{A}$-equivalent Newton-Puiseux parametrization (see Proposition \ref{ordemdasordens}). In Proposition \ref{caracdiffeos}, we characterize the change of coordinates and parameter that preserve this particular parametrization. We define the subgroups $\tilde{\mathcal{A}}_1$ and $\mathcal{H}$ of $\mathcal{A}$ (Definition \ref{definicaoa1tileh}) so that the aforementioned analytic diffeomorphisms can be decomposed into a composition of elements from $\mathcal{H}$ and $\tilde{\mathcal{A}}_1$. In Proposition \ref{acaodosdiffeo}, we explicitly describe the action of an element of $\tilde{\mathcal{A}}_1$, satisfying Proposition \ref{caracdiffeos}, on Newton-Puiseux parametrizations. These results allow us to establish some criteria for eliminating parameters. In Theorem \ref{limpezageral}, we eliminate some parameters using the value semigroup associated to the curve that extend known results for plane branches presented in \cite{Ebey}, \cite{delorme}, \cite{Zariski-book} and \cite{HH2}. As a consequence, we obtain a truncation order of the parametrized curve $\phi$, that is, we obtain a polynomial map-germ $\mathcal{A}$-equivalent to $\phi$.

In \cite{Zariski}, Zariski introduces an analytic invariant for plane branches that can be determined directly from a parametrization, as given in Theorem \ref{limpezageral}. In Theorem \ref{lambdazariski}, we extend this result to any branch in $\mathbb{C}^{n+1}$ under a hypothesis which is naturally fulfilled for the plane case, but not for $n>1$, as we show in Example \ref{exemplo-lambda}.

The presented elimination criteria allow us to establish a connection between an eliminable term of a parametrized curve in $\mathbb{C}^{n+1}$ and an eliminable term of a particular parametrized curve in $\mathbb{C}^{n+2}$ (Theorem \ref{eliminadeumaiffdaoutra}). Based on this result, we can use normal forms for curves in a $d$-dimensional complex space to obtain normal forms in a $(d+1)$-dimensional space.

This paper contributes to the classification problem and description of normal forms for curves in $\mathbb{C}^{n+1}$, as illustrated in Example \ref{ex-3}. Furthermore, we highlight that these results have been applied to construct normal forms for parametric curves of multiplicity 4 in $\mathbb{C}^3$ and $\mathbb{C}^4$ (as presented in \cite{HRoR}).

\section{Numerical semigroups}
\label{Sec-Semigroups}
The theory of numerical semigroups is a useful tool for classification problems concerning analytic curves. In particular, it plays a fundamental role in the topological, bi-Lipschitz, analytical classification of irreducible plane curves, as can be seen, for instance, in \cite{Zariski-book}, \cite{fernandes} and \cite{HH2}.

In what follows, we present some general properties of numerical semigroups that will be explored throughout this text. For more details see \cite{Hefez} and \cite{semigrupo}.

A subset $S$ of $\N$ is called a {\it semigroup} if $0 \in S$ and it is closed under addition. Given any semigroup $S \subseteq \N$, there exists a unique set $\lbrace v_0, v_1, \ldots, v_g \rbrace \subset S$ such that
$$S = \langle v_0, v_1, \ldots, v_g \rangle:=\left\{\sum_{i=0}^{g}\alpha_iv_i;\ \alpha_i\in\mathbb{N}\right\},$$ with $\min(S\setminus\{0\})=v_0 < v_1 < \cdots < v_g$ and $v_i=\min(S\setminus\langle v_0,\ldots ,v_{i-1}\rangle)$ for $1\leq i\leq g$. In  particular, $v_i \not\equiv v_j \mod v_0$ for $i \neq j$ (see \cite[Proposition 6.1]{Hefez}) and $\lbrace v_0, v_1, \ldots, v_g \rbrace$ is contained in any generating set of $S$.
The set $\lbrace v_0, v_1, \ldots, v_g \rbrace$ is called the {\it minimal generating set} of $S$ and $v_0$ is called {\it multiplicity} of $S$.

The elements $l \in \N \setminus S$ are called gaps of $S$. 
Notice that the number of gaps of a semigroup can be finite or not. We have that a semigroup $S = \langle v_0, v_1, \ldots, v_g \rangle$  has a finite number of gaps if and only if $gcd(v_0, v_1, \ldots, v_g) = 1$ or, equivalently, $S$ contains two consecutive elements (see \cite[Proposition 6.2]{Hefez}). In this case, we call $S$ a {\it numerical semigroup}.

In particular, if $S$ is a numerical semigroup, then there exists an element $c \in S$ such that $c - 1 \notin S$ and every natural number $z$ satisfying $z \geq c$ belongs to $S$. The element $c$ is called the {\it conductor} of the semigroup $S$.

In what follows, any semigroup $S$ is considered a numerical semigroup.

The concept of Ap\'ery set of a semigroup $S$ with respect to an element $s\in S$ is a useful tool for describing the elements of $S$ and their gaps.
In this work, we consider the Ap\'ery set of $S$ with respect to its multiplicity.

\begin{definition}\label{definicaoapery}
Let $S \subseteq \N$ be a numerical semigroup with multiplicity $v_0 = \min (S \setminus \lbrace 0 \rbrace)$. 
Let $A : = \lbrace a_0, a_1, \ldots, a_{v_0 - 1} \rbrace\subset S$, where  
$$a_0=0\ \ \mbox{and}\ \ a_i = \min \left\lbrace S \setminus \displaystyle \bigcup _{j = 0}^{i - 1}(a_j + v_0\N) \right\rbrace\ \ \mbox{ for every}\ \ 1\leq i<v_0.$$
The set $A$ is called the Ap\'ery set of $S$ (with respect to $v_0$).
\end{definition}

The elements of the Ap\'ery set of a semigroup $S$ are such that $a_0=0 < a_1< \ldots < a_{v_0-1}$ and satisfy the following properties:

\begin{proposition}(\cite[page 92]{Hefez})\label{propriedadesapery}
Let $S$ be a semigroup with conductor $c$ and $A = \lbrace a_0, a_1, \ldots, a_{v_0 - 1} \rbrace$ its Ap\'ery set. The following properties are satisfied:
\begin{enumerate}
    \item $a_i \not\equiv a_j\mod v_0$ for every $i \neq j$;
    \item $a_i = \min \left( S\cap \{a_i+zv_0;\ z\in\mathbb{Z}\} \right)$;
    \item $S = \displaystyle \bigcup _{i = 0}^{v_0 - 1} (a_i + v_0\N)$;
    \item $c = a_{v_0 - 1} - v_0 + 1$.
\end{enumerate}
\end{proposition}

By the previous proposition, it follows that for any $z\in\mathbb{Z}$ there exists a unique $a_i\in A$ and $k\in\mathbb{Z}$ such that $z=a_i+kv_0$. Moreover, $z\in S$ if and only if $k\geq 0$. This property allows us to describe every element in $S$ and any gap of $S$.

\begin{remark}\label{remarkminimalcontidonoapery}
Let $S$ be a semigroup with minimal generating set $\lbrace v_0, v_1, \ldots, v_g \rbrace$ and Ap\'ery set $A$. By Proposition \ref{propriedadesapery}, item 3, we have that $(A \setminus \lbrace 0 \rbrace) \cup \lbrace v_0 \rbrace$ is a set of generators for $S$, consequently $\{v_1,\ldots ,v_g\}\subseteq A\setminus\{0\}$ and $2\leq g<v_0$.
Moreover, $a_1 = v_1 = \min (A \setminus \lbrace 0 \rbrace) = \min (S \setminus \langle v_0 \rangle)$.
\end{remark}

For an arbitrary numerical semigroup $S$, it is not immediate to obtain the Ap\'ery set and the conductor $c$ of $S$ by means its minimal generating set. For the minimum or the maximum value of $g$, that is, for $g=2$ or $g=v_0-1$, the Ap\'ery set can be easily described (see \cite{semigrupo}). More precisely, if $S=\langle v_0,v_1\rangle$ then $A=\{kv_1;\ 0\leq k<v_0\}$ and $c=(v_0-1)v_1-v_0+1$. For $S=\langle v_0=n,\ldots ,v_{n-1}\rangle$ we get $A=\{0,v_1,\ldots ,v_{n-1}\}$ and $c=v_{n-1}-n+1$. 

\vspace{0.2cm}
In the next proposition, under certain hypothesis, we describe the Ap\'ery set of a numerical semigroup $\Gamma_1$ in terms of the Ap\'ery set of a numerical semigroup $\Gamma_0\subset\Gamma_1$.

\begin{proposition}\label{diminuircond}
Let $\Gamma_0 = \langle v_0, \ldots, v_g \rangle$ be a numerical semigroup and $\lbrace 0, a_1, \ldots, a_{v_0 - 1} \rbrace$ its Ap\'ery set. If $v_{g + 1} \in \mathbb{N} \setminus \Gamma_0$ satisfies $v_{g + 1} \equiv a_{v_0 - 1}\mod v_0$ and $\max \lbrace v_g, a_{v_0 - 2} \rbrace < v_{g+ 1} < a_{v_0 - 1}$, then the Ap\'ery set of $\Gamma_1=\langle v_0,\ldots ,v_g,v_{g+1}\rangle$ is $\lbrace 0, a_1, \ldots, a_{v_0 - 2}, v_{g + 1} \rbrace$ and the conductor of $\Gamma_1$ is $v_{g + 1} - v_0 + 1$.
\end{proposition}
\Dem Since $v_{g + 1}>v_g$ is a gap of $\Gamma_0$ and $\Gamma_0$ is a numerical semigroup, we have that $\Gamma_1=\langle v_0,\ldots ,v_g,v_{g+1}\rangle$ is also a numerical semigroup. By item 2 of Proposition \ref{propriedadesapery},  the conditions $v_{g + 1} < a_{v_0 - 1}$ and $v_{g + 1} \equiv a_{v_0 - 1}\mod v_0$ imply that $v_{g + 1}$ is an element of the Ap\'ery set of $\Gamma_1$.
Moreover, since $a_{v_0 - 2} < v_{g + 1}$, it follows that $a_i$ belongs to the Ap\'ery set of $\Gamma_1$ for every $1\leq i\leq v_{0}-2$, that is, the Ap\'ery set of $\Gamma_1$ is $\lbrace 0, a_1, \ldots, a_{v_0 - 2}, v_{g + 1} \rbrace$. In particular, by item 4 of Proposition \ref{propriedadesapery}, we conclude that the conductor of $\Gamma_1$ is $v_{g + 1} - v_0 +1$.
\cqd

\section{Parametrized Curves in $\mathbb{C}^{n + 1}$}

\label{Sec-Param-Curves}

Denote by $\C \lbrace X, \underline{Y} \rbrace:=\mathbb{C}\{X,Y_1, \ldots, Y_n\}$ the absolutely convergent power series ring in a neighborhood of the origin in $\C ^{n + 1}$ and  $\mathcal{M}_{n + 1} = \langle X, \underline{Y} \rangle : = \langle X, Y_1, \ldots, Y_n \rangle$ its maximal ideal.
Let $\langle f_1,\ldots ,f_r\rangle\in\mathbb{C}\{X,\underline{Y}\}$ be the ideal generated by $f_1, \ldots, f_r \in \mathcal{M}_{n+1}$ such that the local ring $\mathcal{O} = \frac{\C \lbrace X, \underline{Y} \rbrace}{\langle f_1, \ldots, f_r \rangle}\cong\mathbb{C}\{x,y_1,\ldots ,y_n\}$ has Krull dimension one, where $x = X + \langle f_1, \ldots , f_r \rangle$ and $y_i = Y_i + \langle f_1, \ldots , f_r \rangle$, for every $i \in \lbrace 1, \ldots, n \rbrace$.

The germ of the analytic curve in $\C ^{n + 1}$ determined by $f_1, \ldots, f_r$ is the germ of the set $$\mathcal{C} = \lbrace (\alpha _0, \underline{\alpha}) = (\alpha _0, \alpha _1, \ldots, \alpha _n) \in \C ^{n + 1} \, ; \ f_1(\alpha _0, \underline{\alpha}) = \cdots = f_r(\alpha _0, \underline{\alpha}) = 0 \rbrace.$$

We say that the curve determined by $f_1,\ldots ,f_r\in\mathbb{C}\{X,\underline{Y}\}$ is {\it reduced} if the ideal $\langle f_1,\ldots ,f_r\rangle$ is radical. In what follows, every curve is considered reduced.

Two reduced curves $\mathcal{C}_1$ and $\mathcal{C}_2$ in $\C ^{n + 1}$ are called {\it analytically equivalent} if there exists an analytic diffeomorphism $H : U \to V$ such that $H (\mathcal{C}_1 \cap U) = \mathcal{C}_2 \cap V$, where $U$ and $V$ are neighborhoods of the origin in $\C ^{n + 1}$. Equivalently, any such analytic diffeomorphism $H$ induces (and is induced by) an isomorphism $H^*:\mathcal{O}_2\rightarrow\mathcal{O}_1$ between the local rings $\mathcal{O}_2$ and $\mathcal{O}_1$ of $\mathcal{C}_2$ and $\mathcal{C}_1$, respectively. Thus, two reduced curves are analytically equivalent if and only if their local rings are isomorphic as $\mathbb{C}$-algebras (see \cite[Section 8.2]{BK}).

In what follows, we consider the analytic curve $\mathcal{C}$ to be irreducible, which is equivalent to assuming that $\langle f_1,\ldots ,f_r\rangle$ is a prime ideal and consequently that the local ring $\mathcal{O}$ is a domain.

Let $\mathcal{K}$ be the field of fractions of $\mathcal{O}$ and $\overline{\mathcal{O}}$ its integral closure in $\mathcal{K}$. Since $\mathcal{O}$ has Krull dimension one, it follows that $\overline{\mathcal{O}} \cong \C \lbrace t \rbrace$ and $\mathcal{K}\cong\mathbb{C}((t))$, where $\mathbb{C}((t))$ denotes the field of formal Laurent series with coefficients in $\mathbb{C}$. 

The inclusion $\mathcal{O} \subset \overline{\mathcal{O}}$ induces the canonical ring monomorphism 
\begin{equation}\label{monovarphi}
\begin{array}{rccc}
\varphi : & \mathcal{O} & \longrightarrow & \overline{\mathcal{O}} \cong \C \lbrace t \rbrace \\
& g(x, y_1, \ldots , y_n) & \longmapsto & g(p_0(t), p_1(t), \ldots , p_n(t)),
\end{array}
\end{equation}
where $\varphi (x) = p_0(t)$ and $\varphi (y_i) = p_i(t)$ for every $i \in \lbrace 1, \ldots , n \rbrace$. Moreover, we have the canonical valuation on $\mathcal{K}$:
\begin{equation}\label{valuation}
\begin{array}{cccc}
v : & \mathcal{K} & \rightarrow & \overline{\mathbb{Z}}:=\mathbb{Z}\cup\{\infty\} \\
& \frac{g}{h} & \mapsto & ord_t(\varphi(g))-ord_t(\varphi(h))
\end{array}
\end{equation}
for $g,h\in\mathcal{O}$, where $v(0)=\infty$.

Notice that  $0=\varphi (f_i(x, y_1, \ldots, y_n)) = f_i(p_0(t), p_1(t), \ldots, p_n(t))$ for every $1\leq i\leq r$; that is, $(p_0(t), p_1(t), \ldots, p_n(t))$ is a parametrization of $\mathcal{C}$ in a neighborhood of the origin. In what follows, we identify such a parametrization with the smooth map-germ $\phi : (\C, 0) \rightarrow (\C ^{n + 1}, 0)$ given by $$\phi (t) = (p_0(t), p_1(t), \ldots, p_n(t)).$$

\begin{remark}\label{primitive} It follows by (\ref{monovarphi}) that $\mathcal{O} \cong\mathbb{C}\{p_0(t),p_1(t),\ldots ,p_n(t)\}$. Since its field of fractions $\mathcal{K}$ is isomorphic to $\mathbb{C}((t))$, we have that $t=\frac{g}{h}\in\mathcal{K}$ for some $g,h\in\mathcal{O}$. In particular, all exponents in $p_0(t), p_1(t), \ldots, p_n(t)$ do not admit a common divisor distinct to $1$, or equivalently, $\phi(t)$ cannot be reparametrized by a power of a new variable. In this case, we say that $\phi(t)$ is \emph{primitive}.
\end{remark}

\begin{example} 
A parametrization of the form $\phi(t) = (t^{12}-t^{18}, t^{8} + 5t^{18}, t^{14} - \sum_{i\geq 2}t^{10i})$ is not primitive since, by setting $t^2=u$, we obtain a parametrization $\psi(u)=(u^{6} -u^{9}, u^{4}+5u^9, u^{7} -\sum_{i\geq 2}u^{5i})$ satisfying the previous remark.\end{example}

Recall that two curves are analytically equivalent if and only if their local rings are isomorphic. Since the local ring of an irreducible curve is isomorphic to a subalgebra of $\mathbb{C}\{t\}\cong\overline{\mathcal{O}}$, it follows that the analytic classes of irreducible curves coincide with the isomorphism classes of such subalgebras. On the other hand, two subalgebras $\mathcal{O}_1=\mathbb{C}\{p_0(t), p_1(t), \ldots, p_n(t)\}$ and 
$\mathcal{O}_2=\mathbb{C}\{q_0(t), q_1(t), \ldots, q_n(t)\}$ of $\mathbb{C}\{t\}$
are isomorphic if and only if there exists an automorphism $\rho^*(t)$ of $\mathbb{C}\{t\}$ (a change of parameter) and an automorphism  $\sigma^*$ of $\mathbb{C}\{x,y_1,\ldots ,y_n\}$ (a change of generators or coordinates) such that $q_i(t)=\sigma^*(p_i(\rho^*(t)))$ (see \cite{Ebey}). This correspondence induces the following equivalence relation.

We denote by ${\rm Diff} (\C ^l, 0)$ the group of germs of analytic diffeomorphisms of $(\C ^l, 0)$. Consider the Mather group $\mathcal{A} ={\rm Diff}(\C,0)\times
{\rm Diff}(\C ^{n + 1},0)$. Two map-germs $\phi _1, \phi _2 : (\C,0) \to (\C ^{n + 1},0)$ are said to be $\mathcal{A}$-{\it equivalent}, denoted by $\phi _1 \underset{\mathcal{A}}{\sim} \phi _2$, if there exists $(\rho, \sigma) \in \mathcal{A}$ such that $\phi _2 = \sigma \circ
\phi _1 \circ \rho^{-1}$. This condition is illustrated by the following commutative diagram:
$$  \begin{array}{lcr}
  (\C,0) & \stackrel{\phi _1}{\longrightarrow} & (\C ^{n + 1},0) \vspace{0.2cm} \\
  \rho \ \downarrow & \circlearrowleft & \downarrow \ \sigma  \vspace{0.2 cm}\\
 (\C,0)  & \stackrel{\phi _2}{\longrightarrow} &
 (\C ^{n + 1},0)\end{array}$$

In this framework, the study of analytical equivalence between curves reduces to the study of $\mathcal{A}$-equivalence between their parametrizations.

The embedding dimension of the curve, denoted by $e(\O)$, is the dimension of the Zariski tangent space of its local ring, that is, $e(\O)=\dim_{\mathbb{C}}\frac{\mathcal{M}_{\mathcal{O}}}{\mathcal{M}^2_{\mathcal{O}}}$, where $\mathcal{M}_{\mathcal{O}}$ denotes the maximal ideal of $\O$. Thus $e(\mathcal{O})$ is the minimum number of generators for  $\mathcal{M}_{\mathcal{O}}$ (see \cite{Ebey}).

\begin{definition}We say that a parametrization $\phi (t) = (p_0(t), p_1(t), \ldots, p_n(t))$ is \emph{degenerate} if $\phi(t)$ is $\mathcal{A}$-equivalent to $\left( q_0(t), q_1(t), \ldots , q_n(t) \right),
$ with $q_i(t) = 0$ for some $0 \leq i \leq n$. Otherwise, the parametrized curve is called non-degenerate and in this case we say that the \emph{embedding dimension} of the curve parametrized by $\phi(t)$ is $n + 1$.
\end{definition}

\begin{example} 
   The parametrization $\phi _1(t) = (t^{12}, t^{30} + \sum_{i\geq 4}t^{10i}, t^{30} - t^{32}+5t^{33},2t^{32}-10t^{33}+\sum_{i\geq 4}2t^{10i})$ is degenerate. Indeed, taking $\rho(t)=t$ and $\sigma(x,y_1,y_2,y_3)=(x,y_1,y_2,y_3-2(y_1-y_2))$ we obtain $(\sigma\circ\phi_1\circ\rho^{-1})(t)=(t^{12}, t^{30} + \sum_{i\geq 4}t^{10i}, t^{30} - t^{32}+5t^{33},0)$. 
   
   On the other hand, we have that $\phi_2(t)=(t^{6} -t^{9}, t^{4}+5t^9, t^{7} -\sum_{i\geq 2}t^{5i})$ is a non-degenerate parametrization, since $\{t^{6} -t^{9},t^{4}+5t^9, t^{7} -\sum_{i\geq 2}t^{5i}\}\subset\mathcal{M}_{\mathcal{O}}\setminus\mathcal{M}^2_{\mathcal{O}}$ and the class of these elements gives us a basis for $\frac{\mathcal{M}_{\mathcal{O}}}{\mathcal{M}^2_{\mathcal{O}}}$.
\end{example}

Throughout this paper, we consider only irreducible curves in $\mathbb{C}^{n+1}$ with primitive non-degenerate parametrizations.

Let $\mathcal{C} \subset \mathbb{C}^{n+1}$ be a curve parametrized by $\phi : (\C, 0) \rightarrow (\C ^{n + 1}, 0)$. The monomorphism $\varphi$, defined in \eqref{monovarphi}, is such that $\varphi (g) = g(\phi (t))$ for every $g \in \mathcal{O}$. In this way, considering the natural valuation given in (\ref{valuation}) we define the set
\begin{equation}
 \label{Def-Semigroup}  
 \Gamma (\phi) = \lbrace v(g) = ord_t(g(\phi (t))) \, ; \ g \in \mathbb{C}\{X,Y_1,\ldots ,Y_n\} \ \mbox{and} \ g(\phi(t)) \neq 0 \rbrace .
\end{equation}
It follows that $\Gamma (\phi)$ is an additive semigroup, called the {\it values semigroup} associated with $\mathcal{C}$ or $\phi$. Where the context of the parametrization $\phi$ is clear, we shall simply denote the semigroup by $\Gamma$.

\begin{remark}\label{seriescomordemmaiorqueocond}
By Remark \ref{primitive}, there exist $g,h\in\mathcal{O}$ such that $\frac{g}{h}=t\in \mathcal{K}$. Then, $v(g)=v(h)+1$, that is, $v(g)$ and $v(h)$ are two consecutive elements in $\Gamma(\phi)$. Consequently, $\Gamma(\phi)$ has a conductor $c$ (see Section \ref{Sec-Semigroups}). Moreover, the conductor ideal $(\mathcal{O}:\overline{\mathcal{O}}):=\{h\in\mathcal{O};\ h\overline{\mathcal{O}}\subseteq\mathcal{O}\}$ of $\mathcal{O}$ is isomorphic to $t^c\mathbb{C}\{t\}$ (see \cite[Theorem 1]{Ebey}). In particular, any element in $\overline{\mathcal{O}}\cong \mathbb{C}\{t\}$ with order greater than or equal to $c$ belongs to $\mathcal{O}$.
\end{remark}

The values semigroup is invariant under $\A$-equivalence, that is, if $\phi _1 \underset{\mathcal{A}}{\sim} \phi _2$, then $\Gamma(\phi_1)=\Gamma(\phi_2)$ (see \cite{GibsonHobbs}). Given a parametrization $\phi$, the minimal generating set of its semigroup $\Gamma$ can be determined via the algorithm described in \cite{basestandard}.
 
The {\it delta invariant} of $\phi$ is defined by 
$\delta _{\phi} = \dim _{\C} \frac{\overline{\mathcal{O}}}{\mathcal{O}} 
$. Since $\dim_{\mathbb{C}}\frac{\overline{\mathcal{O}}}{\mathcal{O}}=\sharp \left (v(\overline{\mathcal{O}})\setminus v(\mathcal{O})\right )=\sharp\left ( \mathbb{N}\setminus\Gamma\right )$ (see \cite{Kunz}), we obtain that $\delta_{\phi}$ is the number of gaps of the semigroup $\Gamma$  and, consequently, $\delta_{\phi}<\infty$.

Given a map-germ $f : (\C , 0) \to (\C ^{n+1}, 0)$, we denote by $j^kf$ the Taylor polynomial of $f$ at the origin of degree $k$, called the $k$-jet of $f$. We say that $f$ is $\mathcal{A}$-{\it finitely determined} if there exists $k$ such that for every $g : (\C , 0) \to (\C ^{n+1}, 0)$ with $j^k g = j^k f$, then $f \underset{\mathcal{A}}{\sim} g$. In particular, if $f$ is $\mathcal{A}$-finitely determined, then $f \underset{\mathcal{A}}{\sim} j^k f$.

Since $\delta_{\phi}$ is finite, according to \cite[Lemma 3.1,Proposition 4.1]{GibsonHobbs} or \cite[Section 2]{RodriguesRuas}, we have that $\phi$ is $\mathcal{A}$-finitely determined. In Corollary \ref{corolariolimpezageral}, we present an upper bound for $k\in\mathbb{N}$ such that $\phi \underset{\mathcal{A}}{\sim} j^k \phi$ using the semigroup $\Gamma$.

\section{Newton-Puiseux parametrizations}

A solution to the analytic classification problem for irreducible plane curves was presented in  \cite{HH2, handbook}, based on parametrizations expressed in a specific form. In this section, we show that any parametrization of an irreducible non-degenerate curve in $(\mathbb{C}^{n+1},0)$ admits a particular form, namely Newton-Puiseux parametrizations, which share properties with the plane case. This approach allows us to explicitly eliminate terms via $\mathcal{A}$-equivalence by constructing the corresponding changes of parameters and coordinates. 

Let $\overline{\phi}(t)=(q_0(t),q_1(t),\ldots ,q_n(t))$ be a non-degenerate parametrization of an irreducible curve in $(\mathbb{C}^{n+1},0)$. Up to a permutation of the components of $\overline{\phi}$ (which corresponds to the action of an element in $\mathcal{A}$), we may assume that $v_0 = ord _t(q_0(t)) = \min \lbrace ord _t (q_i(t)) \, ; \ 0 \leq i \leq n \rbrace=\min(\Gamma(\overline{\phi})\setminus\{0\})$. Consequently, there exists a unit $u(t) \in \C \lbrace t \rbrace$ such that $q_0(t) = t^{v_0} \cdot u(t).$
By	taking $(\rho , \sigma) \in \mathcal{A}$ defined by $
	\rho (t) = t \cdot u(t)^{\frac{1}{v_0}}$ and $\sigma (X, \underline{Y}) = (X, \underline{Y}),$
	we obtain $q_0(t) = \rho(t)^{v_0}$ and 
    $$
    \phi (t) = 
	(\sigma \circ \overline{\phi} \circ \rho ^{-1})(t) =  \displaystyle \left( t^{v_0}, p_1(t), \ldots, p_n(t) \right),
    $$
	where $p_i(t) = q_i(\rho ^{-1}(t))$, for $i \in \lbrace 1, \ldots, n \rbrace$ with $ord _t(p_i(t)) \geq v_0$. The integer $v_0$ is called the {\it multiplicity} of the curve $\mathcal{C}$ (or of $\phi$). 

Similarly to the values semigroup $\Gamma$ of $\phi(t)=(p_0(t),p_1(t),\ldots ,p_n(t))$, we define 
$$
\Gamma _i  = \left\lbrace ord _t(g(p_0(t), \ldots, p_{i - 1}(t), 0, \ldots, 0)) \, ; \ g \in \mathcal{O} \ \mbox{and} \ g(p_0(t), \ldots, p_{i - 1}(t), 0, \ldots, 0) \neq 0 \right\rbrace $$
\begin{equation}\label{gammai} \hspace{-3.2cm}= \{v(h);\ h\in\mathbb{C}\{x,y_1,\ldots ,y_{i-1}\}\ \mbox{and}\  h(p_0(t), \ldots, p_{i - 1}(t)) \neq 0\}.\end{equation}
Notice that $\Gamma_i$ is a semigroup but it does not  necessarily have a conductor, since the exponents of $p_0(t), \ldots, p_{i - 1}(t)$ may have a common divisor greater than $1$. Where appropriate, we use the notation $\Gamma_i(\phi)$. While preparing the final version of this manuscript, we became aware of a result in \cite[Section 5.2]{mont} similar to the next proposition.

\begin{proposition}\label{ordemdasordens}
Any non-degenerate parametrization of an irreducible curve $\mathcal{C}$ in $\mathbb{C}^{n+1}$ with semigroup $\Gamma= \langle v_0, v_1, \ldots, v_g \rangle$  
is $\mathcal{A}$-equivalent to
\begin{equation}\label{parametgeral2}
\phi(t)=\left( t^{v_0}, t^{w_1} + \sum _{j > w_1} a_{1j}t^j, \ldots, t^{w_n} + \sum _{j > w_n}a_{nj}t^j \right),
\end{equation}
where $\min(\Gamma\setminus\{0\})=v_0 < w_1 < w_2 < \cdots < w_n$ and $w_i = \min (\Gamma \setminus \Gamma _i)$ for every $i \in \lbrace 1, \ldots, n \rbrace$.
In particular, we have $w_1 = \min (\Gamma \setminus \langle v_0 \rangle) = v_1$.
\end{proposition}
\Dem Let $\phi_0:(\mathbb{C},0) \to (\mathbb{C}^{n+1},0)$ be a Newton-Puiseux parametrization of a curve $\mathcal{C}$, given by
$\phi_0 (t) = \displaystyle \left( t^{v_0}, \sum _{j \geq w_{10}}a^{0}_{1j}t^j, \sum _{j \geq w_{20}}a^{0}_{2j}t^j , \ldots , \sum _{j \geq w_{n0}}a^{0}_{nj}t^j \right)$, where $w_{i0}\geq v_0$ for every $1\leq i\leq n$. By the non-degeneracy condition, there exists $w_{i1} = \min \lbrace j\not\in\langle v_0\rangle \, ; \ a^{0}_{ij} \neq 0 \rbrace$ for each $i \in \lbrace 1, \ldots, n \rbrace$. Taking $(\rho, \sigma) \in \mathcal{A}$ with $\rho(t)=t$ and $$
\sigma (X, \underline{Y}) = \left( X, \frac{1}{a^{0}_{1w_{11}}}\left( Y_1 - \sum _{w_{10} \leq j < w_{11}} a^{0}_{1j}X^{\frac{j}{v_0}} \right), \ldots, \frac{1}{a^{0}_{nw_{n1}}}\left( Y_n - \sum _{w_{n0} \leq j < w_{n1}} a^{0}_{nj} X^{\frac{j}{v_0}} \right) \right),
$$ we obtain
$$\phi_1(t):=\sigma\circ\phi_{0}\circ \rho^{-1} (t) = \displaystyle\left( t^{v_0}, t^{w_{11}} + \sum _{j > w_{11}} a^{1}_{1j}t^j,  t^{w_{21}} + \sum _{j > w_{21}} a^{1}_{2j}t^j,\ldots, t^{w_{n1}} + \sum _{j > w_{n1}} a^{1}_{nj}t^j \right)$$ with $w_{i1}\in\Gamma\setminus \langle v_0\rangle=\Gamma\setminus\Gamma_1(\phi_0)$ for every $i \in \lbrace 1, \ldots, n \rbrace$, where $\Gamma_1(\phi_0)$ was defined in (\ref{gammai}). 

We can assume, up to a permutation of the last $n$ components of $\phi_1(t)$ (which corresponds to the action of an element in $\mathcal{A}$), that 
satisfying
$v_0 < w_{11} \leq w_{21} \leq \cdots \leq w_{n1}$. Setting $w_1 := w_{11}$ and $a^{1}_{1j}:=a_{1j}$, notice that $w_1 = v_1$, since $v_1 = \min (\Gamma \setminus \langle v_0 \rangle) = \min (\Gamma \setminus \Gamma_1(\phi_1))$ and $\Gamma_1(\phi_1)=\Gamma_1(\phi_0)$. 

Now, for each $w_{i1}\in\Gamma_2(\phi_1)$ for $2\leq i\leq n$, there exists $h_i\in\mathbb{C}\{x,y_1\}$ such that $j^{w_{i1}}h_i(\phi_1 (t))=t^{w_{i1}}$. Defining $\rho (t) = t$ and 
$\sigma (X, \underline{Y}) = (X, Y_1,\ldots ,Y_{i-1}, Y_i-h_i,Y_{i+1},\ldots, Y_n)$, the non-degeneracy of the parametrization $\phi_0$ ensures that the $i$-th component of $(\sigma \circ \phi_1 \circ \rho ^{-1})(t)$ is non-null. Consequently, this component has order greater than $w_{i1}$. 
Iterating this process with suitable elements $(\rho, \sigma)\in\mathcal{A}$, we find that every $i$-th component (for $2\leq i\leq n$) of the resulting $\mathcal{A}$-equivalent parametrization to $\phi_1$ is given by $\displaystyle \sum _{j \geq w_{i2}} b_{ij}t^j$, with $b_{iw_{i2}} \neq 0$ and $w_{i2} \in\Gamma \setminus \Gamma _2(\phi_1)$.
By applying $(\rho, \sigma) \in \mathcal{A}$ with $\rho(t)=t$ and $
\sigma (X, \underline{Y}) = \left( X,Y_1, \frac{1}{b_{2w_{22}}}Y_2, \ldots, \frac{1}{b_{nw_{n2}}}Y_n \right),$
followed if necessary by a permutation of the last $(n-1)$ components, we have that $\phi_1$ is $\mathcal{A}$-equivalent to 
$$
\phi_2(t) = \left( t^{v_0}, t^{w_1} + \sum _{j > w_1} a_{1j}t^j, t^{w_{22}} + \sum _{j > w_{22}} a^{2}_{2j}t^j, 
\ldots, t^{w_{n2}} + \sum _{j > w_{n2}} a^{2}_{nj}t^j \right)
$$
with $v_0<w_1<w_{22}\leq \ldots \leq w_{n2}$ and $w_{i2}=\min(\Gamma\setminus\Gamma_2(\phi_1))$ for every $2\leq i\leq n$. 
We denote by $w_2 := w_{22}$, and $a^{2}_{2j}:=a_{2j}$. Note that $\Gamma_k(\phi_2)=\Gamma_k(\phi_1)$ for $k\in\{1,2\}$. 

The result follows by repeating this procedure for the remaining components.
\cqd

\begin{remark}\label{observacaodimensaoegenero}
By the previous result, if $\lbrace v_0, v_1, v_2, \ldots, v_g \rbrace$ is the minimal generating set  for $\Gamma$, then $\lbrace v_0, w_1=v_1, w_2, \ldots, w_n \rbrace \subseteq \lbrace v_0, v_1, v_2, \ldots, v_g \rbrace$. Hence, we have 
$$  n + 1 \leq g + 1 \leq v_0,$$
that is, the embedding dimension of a non-degenerate irreducible curve is less than or equal to its multiplicity (see Remark \ref{remarkminimalcontidonoapery}). 
\end{remark}

Since any non-degenerate parametrized curve is $\mathcal{A}$-equivalent to one described in Proposition \ref{ordemdasordens}, we set
$$\begin{array}{ll}
\mathcal{P}_{(v_0, w_1, \ldots, w_n)} = &  \lbrace \phi (t) = (t^{v_0}, y_1(t), \ldots, y_n(t)) \, ; \ j^{w_i} y_i(t) = t^{w_i}, \mbox{with} \\
& \hspace{1cm} w_i=\min(\Gamma\setminus\Gamma_i)\ \mbox{and}\ v_0=w_0<w_1<\ldots <w_n \rbrace .\end{array}$$

The next result characterizes the elements $(\rho,\sigma)\in\mathcal{A}$ such that $(\sigma\circ\phi\circ\rho^{-1})(t)\in \mathcal{P}_{(v_0, w_1, \ldots, w_n)}$ for any $\phi(t)\in\mathcal{P}_{(v_0, w_1, \ldots, w_n)}$. In other words, we characterize the pairs of diffeomorphisms in $\mathcal{A}$ that preserve the initial form of $\phi$.

We call an element in $\mathcal{P}_{(v_0,w_1,\ldots ,w_n)}$ a Newton-Puiseux parametrization.

\begin{proposition}\label{caracdiffeos}
Let $\phi(t)$ be an element in $\mathcal{P}_{(v_0, w_1, \ldots, w_n)}$. A pair $(\rho, \sigma) \in \mathcal{A}$ satisfies $(\sigma\circ\phi\circ\rho^{-1})(t)\in\mathcal{P}_{(v_0, w_1, \ldots, w_n)}$ if and only if \begin{equation}\label{diffeorho}
\rho (t) = t \cdot \left( \epsilon ^{v_0} + \frac{h_0(\phi (t))}{t^{v_0}} \right)^{\frac{1}{v_0}},
\end{equation}
$\sigma(X,\underline{Y}) =(\sigma_0(X,\underline{Y}),\ldots ,\sigma_n(X,\underline{Y}))$ with
\begin{equation}\label{diffeosigma}
\sigma _0 (X, \underline{Y})  =  \epsilon ^{v_0}X + h_0(X, \underline{Y}) \ \ \mbox{and}\ \
\sigma _k(X, \underline{Y})  = \epsilon ^{w_k} Y_k + h_k (X, \underline{Y})
\end{equation}
where $\epsilon \in \mathbb{C} ^*$ and $h_k \in \langle Y_{k + 1}, \ldots , Y_n \rangle + \mathcal{M}_{n + 1}^2$ is such that $ord _t(h_k(\phi (t))) > w_k$ for $0\leq k\leq n$. 
\end{proposition}
\Dem Let $$\phi(t)=\left( t^{v_0}, t^{w_1} + \sum _{j > w_1} a_{1j}t^j, \ldots, t^{w_n} + \sum _{j > w_n}a_{nj}t^j \right)\in \mathcal{P} _{(v_0, w_1, \ldots, w_n)}.$$

If a pair $(\rho,\sigma)\in\mathcal{A}$ satisfies (\ref{diffeorho}) and (\ref{diffeosigma}), then 
\begin{equation}
\label{sigma-k-phi-rho}(\sigma\circ\phi\circ\rho^{-1})(t)=\left(\ \epsilon^{v_0}\rho^{-1}(t)^{v_0}+h_0(\phi(\rho^{-1}(t))),\ \sigma_1(\phi(\rho^{-1}(t))), \ldots, \ \sigma_n (\phi(\rho^{-1}(t)))\ \right),\end{equation}
where $\sigma_k(\phi(\rho^{-1}(t)))=\epsilon^{w_k}(\rho^{-1}(t)^{w_k}+ \sum _{j > w_k} a_{kj}\rho^{-1}(t)^j)+h_k(\phi(\rho^{-1}(t))$ for $1 \leq k \leq n$.

\vspace{0.2cm}
Since $\rho(t)^{v_0}=\epsilon^{v_0}t^{v_0}+h_0(\phi(t))$ and $j^1(\rho^{-1}(t))=\epsilon^{-1}t$, we obtain $$(\sigma\circ\phi\circ\rho^{-1})(t)=\left( t^{v_0}, t^{w_1} + \sum _{j > w_1} b_{1j}t^j, \ldots, t^{w_n} + \sum _{j > w_n}b_{nj}t^j \right)\in \mathcal{P} _{(v_0, w_1, \ldots, w_n)}.$$

Reciprocally, for any $(\rho,\sigma)\in\mathcal{A}$, we can write $\rho(t)=t\cdot u(t)$ and 
$\sigma (X, \underline{Y}) =(\sigma_0 (X, \underline{Y}),\ldots ,\sigma_n (X, \underline{Y}))$ with 
$$\sigma_k(X,\underline{Y})= \alpha_{k0}X + \sum _{i = 1}^n \alpha_{ki}Y_i +p_k(X,\underline{Y}),
$$
where $u(t)\in\mathbb{C}\{t\}$ is a unit, $p_k(X,\underline{Y})\in\mathcal{M}^2_{n+1}$ and $\alpha_{ki} \in \mathbb{C}$, for $0\leq k,i\leq n$. 
In this way, we get $$(\sigma \circ \phi \circ \rho ^{-1})(t)=
\sigma \left(  (\rho ^{-1}(t))^{v_0}, (\rho ^{-1}(t))^{w_1} + \sum _{j > w_1} a_{1j}(\rho ^{-1}(t))^j, \ldots, (\rho ^{-1}(t))^{w_n} + \sum _{j > w_n}a_{nj}(\rho ^{-1}(t))^j \right) 
$$
where its $k$-th component is given by
\begin{equation}\label{sigmak}
(\sigma_k\circ\phi\circ\rho^{-1})(t)= \alpha_{k0} (\rho ^{-1}(t))^{v_0} + \sum _{i = 1}^n \alpha_{ki} \left( \rho ^{-1}(t)^{w_i} + \sum _{j > w_i} a_{ij}(\rho ^{-1}(t))^j \right) +p_k(\phi(\rho^{-1}(t))).
\end{equation}

Suppose $(\rho,\sigma)\in\mathcal{A}$ such that $(\sigma\circ\phi\circ\rho^{-1})(t)\in\mathcal{P} _{(v_0, w_1, \ldots, w_n)}$. 

Since $\sigma_0(X,\underline{Y})=\alpha_{00}X+h_0(X,\underline{Y})$, where $h _0 (X, \underline{Y}) := \displaystyle \sum _{i = 1}^n \alpha_{0i}Y_i + p_0(X,\underline{Y})\in \langle Y_1, \ldots, Y_n \rangle + \mathcal{M}_{n + 1}^2$, we obtain $(\sigma _0\circ\phi\circ \rho ^{-1})(t)=
\alpha_{00} (\rho ^{-1}(t))^{v_0} + h _0(\phi (\rho ^{-1}(t)))$. Thus, $\alpha _{00} (\rho ^{-1}(t))^{v_0} + h _0(\phi (\rho ^{-1}(t))) = t^{v_0}
$ implies that $$\rho (t) = \left( \alpha _{00}t^{v_0} + h _0(\phi (t)) \right)^{\frac{1}{v_0}} = t \cdot \left( \alpha _{00} + \frac{h _0(\phi (t))}{t^{v_0}} \right)^{\frac{1}{v_0}},
$$
where $\alpha _{00} \neq 0$ because $\rho(t) \in \mbox{Diff}(\mathbb{C} , 0)$. By Proposition \ref{ordemdasordens}, we must have $ord _t(h _0(\phi (t)) > v_0$.
Moreover, taking $\epsilon  \in \mathbb{C} ^*$ such that $\epsilon ^{v_0}=\alpha_{00}$ we obtain
\begin{equation}\label{exp-rho}
\rho (t) = t \cdot \left( \epsilon ^{v_0} + \frac{h _0(\phi (t))}{t^{v_0}} \right)^{\frac{1}{v_0}} \ \mbox{and} \ \sigma _0(X, \underline{Y}) = \epsilon ^{v_0}X + h _0(X, \underline{Y}),
\end{equation} with $h _0 \in \langle Y_1, \ldots, Y_n \rangle + \mathcal{M}_{n + 1}^2$.

In addition, if $(\sigma \circ \phi \circ \rho ^{-1})(t) \in \mathcal{P}_{(v_0, w_1, \ldots, w_n)}$, then for $1\leq k\leq n$ it follows that 
$$j^{w_k}\left((\sigma_k \circ \phi \circ \rho^{-1})(t)\right)= t^{w_k},$$
where $(\sigma_k \circ \phi \circ \rho^{-1})(t)$ is given in (\ref{sigmak}).

By (\ref{exp-rho}), we have $j^1(\rho(t))=\epsilon\cdot t$ and consequently
$ord _t(\rho ^{-1}(t)) = ord _t(\rho (t)) = 1$. Since $v_0 < w_1 < w_2 < \cdots < w_n$ it follows that $\alpha _{ki} =0$ for every $0\leq i<k$ and $\alpha _{kk} \neq 0$, that is, $$
\sigma _{k}(X, \underline{Y}) = \alpha_{kk}Y_k + \sum _{i = k + 1}^n \alpha_{ki}Y_i + p_k(X,\underline{Y}).
$$
 Denoting $\displaystyle h _k (X, \underline{Y}) := \sum _{i = k + 1}^n \alpha_{ki}Y_i +p_k(X,\underline{Y}) \in \langle Y_{k + 1}, \ldots, Y_n \rangle + \mathcal{M}_{n + 1}^2$, we obtain $
\sigma _k(X, \underline{Y}) = \alpha_{kk} Y_k + h _k(X, \underline{Y}).
$
Moreover, from the equality $j^{w_k}\left((\sigma_k \circ \phi \circ \rho^{-1})(t)\right)= t^{w_k}$, or equivalently $j^{w_k}\left((\sigma_k \circ \phi (t)\right)=j^k(\rho(t))^{w_k}$, and using (\ref{exp-rho}), we have 
$$
j^{w_k}\left (\alpha_{kk} \cdot \left( t^{w_k} + \sum _{j > w_k} a_{kj}t^j \right) + h _k (\phi (t)) \right )= j^{w_k}\left (\rho (t)^{w_k}\right )=\epsilon^{w_k}\cdot t^{w_k}.$$
Therefore, we must have $\alpha_{kk} = \epsilon ^{w_k}$, hence
$$
\sigma _k(X, \underline{Y}) = \epsilon ^{w_k} Y_k + h _k (X, \underline{Y}),
$$
where $h _k \in \langle Y_{k + 1}, \ldots, Y_n \rangle + \mathcal{M}_{n + 1}^2$ is such that $ord _t(h _k(\phi (t))) > w_k$.
\cqd 

\begin{definition}\label{definicaoa1tileh}
We denote by $\tilde{\mathcal{A}}_1$ the subgroup of $\mathcal{A}$ consisting of elements $(\rho, \sigma) \in \mathcal{A}$ such that $j^1 \rho (t) = t$ and $$
j^1 \sigma (X, \underline{Y}) = \left(X + \sum _{j = 1}^n\alpha_{1j} Y_j, Y_1 + \sum _{j = 2}^n \alpha_{2j}Y_j, Y_2 + \sum _{j = 3}^n \alpha_{3j}Y_j, \ldots, Y_{n-1}+\alpha_{n-1,j}Y_n, Y_n \right),
$$ meaning that the Jacobian matrix of $\sigma$ at the origin is unipotent upper triangular. In the particular case where $\alpha_{ij} = 0$ for every $i$ and $j$, we obtain the well-known subgroup $\mathcal{A}_1 \subset \mathcal{A}$ of diffeomorphisms whose $1$-jet is the identity.
In addition, we denote by $\mathcal{H}$ the $\mathcal{A}$-subgroup of the homotheties, that is, \begin{equation}\label{grupoH}
\mathcal{H} = \lbrace (\rho, \sigma) \in \mathcal{A} \, ; \, \rho (t) = \alpha t \ \mbox{and} \ \sigma (X, \underline{Y}) = (\alpha _0X, \alpha _1Y_1, \ldots, \alpha _nY_n), \alpha, \alpha _i \in \mathbb{C} ^*, 0 \leq i \leq n \rbrace .
\end{equation}
\end{definition}

Notice that each pair of diffeomorphisms $(\rho , \sigma)$ in Proposition \ref{caracdiffeos} can be expressed as a composition of an element in $\mathcal{H}$ with an element in $\tilde{\mathcal{A}}_1$.

In the next result, we explicitly describe the action of an element of the subgroup $\tilde{\mathcal{A}}_1$ satisfying Proposition \ref{caracdiffeos}.

\begin{proposition}\label{acaodosdiffeo}
Let $\phi (t)=\left( x(t),y_1(t),\ldots ,y_n(t)\right )$ be an element in $\mathcal{P} _{(v_0, w_1, \ldots, w_n)}$, where $x(t)=t^{v_0}$ and $y_i(t)=t^{w_i} + \sum _{j > w_i} a_{ij}t^j$ for each $1\leq i\leq n$. Consider $(\rho , \sigma) \in \tilde{\mathcal{A}}_1$ satisfying Proposition \ref{caracdiffeos}.
Then, $$(\sigma \circ \phi \circ \rho ^{-1})(t) = \phi (t) + \left( 0, \theta _1(\rho ^{-1}(t)), \ldots,  \theta _n(\rho ^{-1}(t)) \right),$$ where $\displaystyle  \theta _i(t) = \displaystyle  h _i(\phi (t)) - \frac{h _0(\phi (t)) y_i ' (t)}{x'(t)} - \sum _{j \geq w_ i}a_{ij}t^j \sum _{l = 2}^{\infty} \left( \begin{array}{c}
\frac{j}{v_0} \\
l
\end{array} \right) \left( \frac{h _0(\phi (t))}{t^{v_0}} \right)^l$. 
Moreover, the initial term of $\theta _i (\rho ^{-1}(t))$ coincides with that of $\theta _i (t)$.
\end{proposition}
\Dem 
Let $(\rho, \sigma) \in \tilde{\mathcal{A}}_1$ satisfying Proposition \ref{caracdiffeos}. Then, $\epsilon = 1$ and we get $$\displaystyle
\rho (t) = t \left(  1 + \frac{h _0(\phi (t))}{t^{v_0}} \right)^{\frac{1}{v_0}}\hspace{0.4cm}\mbox{and} \hspace{0.4cm}
\sigma (X, \underline{Y}) = (\sigma _0(X, \underline{Y}), \ldots, \sigma _n(X, \underline{Y})),$$
where $\sigma _0(X, \underline{Y}) = X + h _0(X, \underline{Y})$ and $\sigma _i(X, \underline{Y}) = Y_i + h _i (X, \underline{Y}),$
with $h _i \in \langle Y_{i + 1}, \ldots, Y_n \rangle + \mathcal{M}_{n + 1}^2$ such that $ord_t(h_i(\phi(t))>w_i$, for every $0\leq i\leq n$.

According to (\ref{sigma-k-phi-rho}), we have $(\sigma \circ \phi \circ \rho ^{-1})(t) = (t^{v_0}, \sigma _1(\phi (\rho ^{-1}(t))), \ldots, \sigma _n(\phi (\rho ^{-1}(t))))$, where
\begin{equation}\label{sigmaiauxiliar}
(\sigma _i \circ\phi\circ \rho ^{-1})(t) = y_i(\rho ^{-1}(t)) + h _i(\phi (\rho ^{-1}(t))),
\end{equation}
for $1\leq i\leq n$.
Notice that for $k \geq 0$ we get $$
\begin{array}{rcl}
\rho (t)^k = \displaystyle t^k \cdot \left( 1 + \frac{h _0(\phi (t))}{t^{v_0}} \right)^{\frac{k}{v_0}} & = & \displaystyle t^k + t^k \cdot \sum _{l = 1}^{\infty} \left( \begin{array}{c}
\frac{k}{v_0} \\
l
\end{array} \right) \left( \frac{h _0(\phi (t))}{t^{v_0}} \right)^l.
\end{array}
$$
Therefore, $$
(\rho ^{-1}(t))^k = t^k - (\rho ^{-1}(t))^k  \sum _{l = 1}^{\infty} \left( \begin{array}{c}
\frac{k}{v_0} \\
l
\end{array} \right) \left( \frac{h _0(\phi (\rho ^{-1}(t)))}{( \rho ^{-1}(t))^{v_0}} \right)^l.
$$
Since $y_i(t) = \displaystyle \sum _{j \geq w_i} a_{ij}t^j$, with $a_{iw_i}=1$ for $1 \leq i \leq n$, we obtain
$$\begin{array}{rcl}
\displaystyle y_i(\rho ^{-1}(t)) & = & \displaystyle \sum _{j \geq w_i} a_{ij} \left( t^j - (\rho ^{-1}(t))^j \sum _{l = 1}^{\infty} \left( \begin{array}{c}
\frac{j}{v_0} \\
l
\end{array} \right) \left( \frac{h _0(\phi (\rho ^{-1}(t)))}{( \rho ^{-1}(t))^{v_0}} \right)^l \right) \vspace{0.2cm} \\
& = & \displaystyle \sum _{j \geq w_i} a_{ij}t^j - \sum _{j \geq w_i} a_{ij}  (\rho ^{-1}(t))^j \sum _{l = 1}^{\infty} \left( \begin{array}{c}
\frac{j}{v_0} \\
l
\end{array} \right) \left( \frac{h _0(\phi (\rho ^{-1}(t)))}{( \rho ^{-1}(t))^{v_0}} \right)^l \vspace{0.2cm} \\
& = & \displaystyle y_i(t) - \sum _{j \geq w_i} a_{ij}  (\rho ^{-1}(t))^j \sum _{l = 1}^{\infty} \left( \begin{array}{c}
\frac{j}{v_0} \\
l
\end{array} \right) \left( \frac{h _0(\phi (\rho ^{-1}(t)))}{( \rho ^{-1}(t))^{v_0}} \right)^l. \vspace{0.2cm} \\
\end{array}
$$
In this way, by equality \eqref{sigmaiauxiliar} we have $$
(\sigma _i\circ\phi \circ\rho ^{-1})(t) = y_i(t) + h _i(\phi (\rho ^{-1}(t))) - \sum _{j \geq w_i} a_{ij}  (\rho ^{-1}(t))^j \sum _{l = 1}^{\infty} \left( \begin{array}{c}
\frac{j}{v_0} \\
l
\end{array} \right) \left( \frac{h _0(\phi (\rho ^{-1}(t)))}{( \rho ^{-1}(t))^{v_0}} \right)^l.
$$
To conclude the result, we highlight the term in the sum $$\displaystyle \sum _{j \geq w_i} a_{ij} (\rho ^{-1}(t))^j \sum _{l = 1}^{\infty} \left( \begin{array}{c}
\frac{j}{v_0} \\
l
\end{array} \right) \left( \frac{h _0(\phi (\rho ^{-1}(t)))}{( \rho ^{-1}(t))^{v_0}} \right)^l$$
given by $l = 1$, that is, 
$$
\begin{array}{rcl}
\displaystyle \sum _{j \geq w_i} a_{ij}(\rho ^{-1}(t))^j \cdot \frac{j}{v_0} \cdot \left( \frac{h _0(\phi(\rho ^{-1}(t)))}{(\rho ^{-1}(t))^{v_0}} \right) & = & \displaystyle \sum _{j \geq w_i} \frac{j a_{ij}(\rho^{-1}(t))^j}{v_0(\rho ^{-1}(t))^{v_0}} \cdot h _0(\phi (\rho^{-1}(t))) \vspace{0.2cm} \\
& = & \displaystyle \sum _{j \geq w_i} \frac{j a_{ij}(\rho^{-1}(t))^{j - 1}}{v_0(\rho ^{-1}(t))^{v_0 - 1}} \cdot  h _0(\phi (\rho^{-1}(t))) \vspace{0.2cm} \\
& = & h _0(\phi(\rho^{-1}(t)) \displaystyle \frac{y_i'(\rho^{-1}(t))}{x'(\rho^{-1}(t))}.
\end{array}
$$
Thus, $$
\begin{array}{rcl}
(\sigma _i\circ\phi\circ \rho ^{-1})(t) & = & \displaystyle y_i(t) + h _i(\phi (\rho ^{-1}(t))) - h _0(\phi (\rho ^{-1}(t))) \frac{y_i'(\rho ^{-1}(t))}{x'(\rho ^{-1}(t))} \vspace{0.2cm} \\
& & - \displaystyle \sum _{j \geq w_i} a_{ij} (\rho ^{-1}(t))^j \sum _{l = 2}^{\infty} \left( \begin{array}{c}
\frac{j}{v_0} \\
l
\end{array} \right) \left( \frac{h _0(\phi (\rho ^{-1}(t)))}{( \rho ^{-1}(t))^{v_0}} \right)^l \vspace{0.2cm} \\
& = & y_i(t) + \theta _i(\rho ^{-1}(t)),
\end{array}
$$
where  $$\theta _i(t) = h _i(\phi (t)) - h _0(\phi (t)) \frac{y_i'(t)}{x'(t)} - \displaystyle \sum _{j \geq w_i} a_{ij} t^j \sum _{l = 2}^{\infty} \left( \begin{array}{c}
\frac{j}{v_0} \\
l
\end{array} \right) \left( \frac{h _0(\phi (t))}{t^{v_0}} \right)^l\ \ \mbox{for}\ 1\leq i\leq n.$$

Therefore, $(\sigma \circ \phi \circ \rho ^{-1})(t) = \phi (t) + \left(0, \theta _1(\rho ^{-1}(t)), \ldots, \theta _n(\rho ^{-1}(t))\right)$.
Moreover, since $j^1(\rho ^{-1}(t)) = t$ we see that $\theta _i(t)$ and $\theta _i(\rho ^{-1}(t))$ have the same initial term.
\cqd

\section{Some elimination criteria}
Based on the results in the previous section, we present some elimination criteria for parameters in order to obtain a simplified parametrization. 

Let $\phi(t)=(t^{v_0},y_1(t),\ldots, y_n(t))\in \mathcal{P}_{(v_0,w_1,\ldots ,w_n)}$ with $y_i(t)=t^{w_i}+\sum_{j>w_i}a_{ij}t^j$ for $1\leq i\leq n$ and
let $(\rho , \sigma) \in \tilde{\mathcal{A}}_1$ satisfying Proposition \ref{caracdiffeos}. According to Proposition \ref{acaodosdiffeo}, we get $(\sigma \circ \phi \circ \rho ^{-1})(t) = \phi (t) + \left( 0, \theta _1(\rho ^{-1}(t)), \ldots,  \theta _n(\rho ^{-1}(t)) \right)$, where 
\begin{equation}\label{theta-i}\displaystyle  \theta _i(t) = \displaystyle  h _i(\phi (t)) - \frac{h _0(\phi (t)) y_i ' (t)}{x'(t)} - \sum _{j \geq w_ i}a_{ij}t^j \sum _{l = 2}^{\infty} \left( \begin{array}{c}
\frac{j}{v_0} \\
l
\end{array} \right) \left( \frac{h _0(\phi (t))}{t^{v_0}} \right)^l.
\end{equation}

In what follows, we consider the $\mathcal{O}$-module $\Omega_i:=\left \{\eta_i-\eta_0\frac{dy_i}{dx};\ \eta_0,\eta_i\in\mathcal{O}\right \}$. Given $\omega=\eta_i-\eta_0\frac{dy_i}{dx}\in\Omega_i$ we define $$v(\omega):=ord_t \left (\eta_i-\eta_0\frac{y'_i(t)}{x'(t)}\right )$$
where $\mathcal{O} \cong\mathbb{C}\{x(t),y_1(t),\ldots ,y_n(t)\}$ with $x(t)=t^{v_0}$.

\begin{remark}\label{remarkeliminacao} By Proposition \ref{acaodosdiffeo}, if $k = ord _t(\theta _i(t))$, then there exist $h_0\in \langle Y_1, \ldots , Y_n \rangle + \mathcal{M}_{n + 1}^2$ and $h_i \in \langle Y_{i + 1}, \ldots , Y_n \rangle + \mathcal{M}_{n + 1}^2$ such that $j^k(\theta_i(\rho^{-1}(t)))=-a_{ik}t^k$. Consequently, $\phi(t)$ is $\mathcal{A}$-equivalent to $\overline{\phi}(t) = (t^{v_0}, \overline{y}_1(t), \ldots, \overline{y}_n(t))$, where 
$j^k (\overline{y}_i(t)) = j^{k-1}(y_i(t))$. In this case, we say that the term of order $k$ in $y_i(t)$ is \emph{eliminable}. Note that the choice of $h_0$ can modify all components of $\phi(t)$ by adding $\theta_i$ to the $i$-th component, but $h_j$ for $j>0$ just change the $j$-th component of $\phi(t)$.

The same analysis is performed if $k=ord _t \left( \theta _i(t) \right) = v \left( h _i(\phi(t)) - h _0(\phi(t)) \frac{y_i'(t)}{x'(t)} \right)$. In this case, the set of orders $$\Lambda_i:=\{v(\omega);\ \omega\in\Omega_i\}$$ is closely related to the eliminable terms in $y_i(t)$. The set $\Lambda_i$ can be computed by the algorithms presented in \cite{basestandard}. The link between coordinate changes and differentials, originally studied by Delorme \cite{delorme} for plane branches, was connected to the tangent space to the $\tilde{\mathcal{A}}_1$-orbit for a plane curve in \cite{HH2} and \cite{RodriguesHernandes}. Furthermore, in \cite[Proposition 2]{HRR}, elements of $\Omega_i$ are related to elements in the extended tangent space to the $\mathcal{A}$-orbit of $\phi:(\mathbb{C},0) \to (\mathbb{C}^{n+1},0)$. 
\end{remark}

In the next result, we identify some eliminable terms of a parametrization by means its semigroup.

\begin{lemma}\label{lemmaeliminasemigrupomaisv1menosv0esemigrupo}
Let $\mathcal{C}$ be a curve with parametrization $\displaystyle \phi (t) =(t^{v_0},y_1(t),\ldots ,y_n(t))\in\mathcal{P}_{(v_0,w_1,\ldots ,w_n)}$ with $y_i(t)= t^{w_i} + \sum _{j > w_i} a_{ij}t^j$ for $1\leq i\leq n$ and semigroup $\Gamma = \langle v_0, v_1, \ldots, v_g \rangle$ where $w_1=v_1$.
\begin{enumerate}
\item If $a_{ik}\neq 0$ for some $k \in \Gamma$ and some $1\leq i\leq n$, then $\phi(t)$ is $\mathcal{A}$-equivalent to $\overline{\phi}(t) = \left( t^{v_0}, y_1(t), \ldots, y_{i-1}(t),  \overline{y}_i(t), y_{i+1}(t),\ldots,  y_n(t) \right)$, with $j^k(\overline{y}_i(t)) = j^{k - 1}(y_i(t))$.

    \item If $a_{ik}\neq 0$ for some $k \in \Gamma + w_i - v_0$ and $1\leq i\leq n$, then $\phi(t)$ is $\mathcal{A}$-equivalent to $\overline{\phi}(t) = (t^{v_0}, \overline{y}_1(t), \ldots, \overline{y}_n(t))$, where $j^k( \overline{y}_i(t)) = j^{k - 1}( y_i(t))$ and $j^{k+w_s-w_i-1}( \overline{y}_s(t)) = j^{k+w_s-w_i-1}( y_s(t))$ for every $s>i$.
\end{enumerate}
\end{lemma}
\Dem 
\begin{enumerate}
    \item If $a_{ik}\neq 0$ for some $k\in\Gamma$, then $k > w_i$ and there exists $h_i\in\langle Y_{i+1},\ldots ,Y_{n}\rangle+\mathcal{M}^2_{n+1}$ such that $v(h_i)=k$ (see (\ref{Def-Semigroup})). In this way, taking $h_s=0$ for every $0\leq s\leq n$ with $s\neq i$ and the associated element $(\rho,\sigma)\in\tilde{\mathcal{A}}_1$ satisfying Proposition \ref{caracdiffeos}, we get $ord_t(\theta_i(t))=ord_t(h_i(\phi(t))=k$ and $\theta_s(t)=0$ for every $1\leq s\leq n$ with $s\neq i$.  It follows, by Remark \ref{remarkeliminacao}, that $\phi(t)$ is $\mathcal{A}$-equivalent to $\overline{\phi}(t) = \left( t^{v_0}, y_1(t), \ldots, y_{i-1}(t),  \overline{y}_i(t), y_{i+1}(t),\ldots,  y_n(t) \right)$ with $j^{k}(\overline{y}_i(t)) = j^{k - 1}(y_i(t))$.

\item If $a_{ik}\neq 0$ for some $k\in\Gamma+w_i-v_0$, then $k>w_i$ and $v_0<k-w_i+v_0$. Thus, there exists $h_0\in\langle Y_{1},\ldots ,Y_{n}\rangle+\mathcal{M}^2_{n+1}$ such that $k=v(h_0)+w_i-v_0$, or equivalently, $k-w_i+v_0=v(h_0)\in\Gamma$. Thus, taking $(\rho,\sigma)\in\tilde{\mathcal{A}}_1$ from Proposition \ref{caracdiffeos}, where $h_s=0$ for every $1\leq s\leq n$, we obtain $$ord_t(\theta_s(t))=ord_t\left (h_0(\phi(t))\frac{y'_s(t)}{x'(t)}\right )=v(h_0)+w_s - v_0=k+w_s-w_i,$$ since $v(h_0)>v_0$ and the order of the sum in (\ref{theta-i}) is at least $v(h_0)+2(w_s-v_0)=k+w_s-w_i+v(h_0)-v_0>k+w_s-w_i$. In this way, by Remark \ref{remarkeliminacao}, we conclude that $\phi(t)$ is $\mathcal{A}$-equivalent to $\overline{\phi}(t) = \left( t^{v_0}, \overline{y}_1(t), \ldots, \overline{y}_n(t) \right)$ with $j^{k}(\overline{y}_i(t)) = j^{k - 1}(y_i(t))$ and $j^{k+w_s-w_i-1}( \overline{y}_s(t)) = j^{k+w_s-w_i-1}( y_s(t))$ for every $s>i$.

\end{enumerate}
\cqd 

A suitable parameter elimination criterion is obtained when we require that the elimination of a term $t^k$ from a component $y_i$ preserves the $(k-1)$-jet of $\phi$. According to item 1 of the previous lemma, the elements $(\rho, \sigma)\in\tilde{\mathcal{A}}_1$ used to eliminate the term $t^k$ (with $k\in\Gamma$ and $k >w_i$) do not affect the other components. In contrast, the change of parameters and coordinates in item 2 can introduce terms of order less than $k$ in the $s$-th component for $s < i$. Thus, to provide a systematic method for eliminating parameters, we apply the previous lemma, and specifically item 2 for the case $i=1$.

As a consequence of Lemma \ref{lemmaeliminasemigrupomaisv1menosv0esemigrupo}, we obtain the following result.

\begin{theorem}\label{limpezageral}
Any non-degenerate parametrization of an irreducible curve $\mathcal{C}$ in $\C ^{n + 1}$ with semigroup $\Gamma = \langle v_0, v_1, \ldots, v_g \rangle$ is $\mathcal{A}$-equivalent to
\begin{equation}\label{curvagerallimpa}
    \left( t^{v_0}, t^{w_1} + \underset{j \notin \Gamma \cup (\Gamma + v_1 - v_0)}{\sum _{j > w_1}} a_{ij}t^j, t^{w_2} + \underset{j \notin \Gamma}{\sum _{j > w_2}} a_{2j}t^j , \cdots , t^{w_n} + \underset{j \notin \Gamma}{\sum _{j > w_n}} a_{nj}t^j \right),
\end{equation}
where $v_0 < w_1 < w_2 < \cdots < w_n$, $w_1 = v_1 = \min (\Gamma \setminus \langle v_0 \rangle)$ and $w_i = \min (\Gamma \setminus \Gamma _{i})$, for $i \in \lbrace 2, \ldots, n \rbrace$ with $\Gamma _i$ given in \eqref{gammai}.
\end{theorem}
\Dem 
By Proposition \ref{ordemdasordens}, any curve with semigroup $\Gamma$ admits a parametrization $\mathcal{A}$-equivalent to
$$
\psi (t) = \left( t^{v_0}, t^{w_1} + \sum _{j>w_1}b_{1j}t^j, \ldots, t^{w_n} + \sum _{j>w_n} b_{nj}t^j \right )\in\mathcal{P}_{(v_0,w_1,\ldots ,w_n)}.
$$

Let $c$ be the conductor of $\Gamma$. According to Lemma \ref{lemmaeliminasemigrupomaisv1menosv0esemigrupo}, we can eliminate every term $b_{1k}t^{k}$ for $w_1<k<c$ with $k\in\Gamma+v_1-v_0$, and every term $b_{ik}t^k$ for $w_i<k<c$ with $k\in\Gamma$, for $1\leq i\leq n$.

Thus, $\psi(t)$ is $\mathcal{A}$-equivalent to the parametrization $\phi(t)=(t^{v_0},y_1(t),\ldots ,y_n(t))$, where
$$\begin{array}{l}
y_1(t) = \displaystyle t^{w_1} + \underset{j \notin \Gamma \cup (\Gamma + v_1 - v_0)}{\sum _{w_1 < j < d_1}} a_{1j}t^j + \sum _{j \geq d_1} a_{1j}t^j\ \ \ \ \mbox{and}\\ \\
y_i(t) = \displaystyle t^{w_i} + \underset{j \notin \Gamma}{\sum _{w_1 < j < d_i}} a_{ij}t^j + \sum _{j \geq d_i} a_{ij}t^j \ \ \mbox{for}\ \ 2\leq i\leq n
\end{array}
$$
where $d_i = \max \lbrace w_i + 1, c \rbrace$ for $1\leq i\leq n$.

Since $k \in \Gamma$ for any $k\geq d_i \geq c$, it follows from Remark \ref{seriescomordemmaiorqueocond} that $$\sum_{j\geq d_i}a_{ij}t^j=t^{d_i}\left(\sum_{j\geq d_i}a_{ij}t^{j-d_i}\right) \in (\mathcal{O}:\overline{\mathcal{O}})\simeq t^c \mathbb{C}\{t\}\ \mbox{for every}\ 1\leq i\leq n.$$ Thus, as $k>w_i$ there exists $h _i \in \langle Y_{i + 1}, \ldots, Y_n \rangle + \mathcal{M}_{n + 1}^2$ such that $
h _i(\phi (t)) = - \sum _{j \geq d_i} a_{ij}t^j.
$ 
Now, let $(\rho,\sigma)\in\tilde{\mathcal{A}}_1$ provided by Proposition \ref{caracdiffeos} with $\rho(t)=t$ (that is $h_0=0$) and $\sigma(X,\underline{Y})=(X,Y_1+h_1(X,\underline{Y}), \ldots, Y_n+ h_n(X,\underline{Y}))$. Consequently, by Proposition \ref{acaodosdiffeo}, we get
$$\begin{array}{ll}
(\sigma\circ\phi\circ\rho^{-1})(t) & = \phi(t)+(0,\theta_1(t),\ldots ,\theta_n(t)) \\
& =\phi(t)+\left ( 0,- \sum _{j \geq d_1} a_{1j}t^j,\ldots , - \sum _{j \geq d_n} a_{nj}t^j\right ) \\
& 
=\left( t^{v_0}, t^{w_1} + \underset{j \notin \Gamma \cup (\Gamma + v_1 - v_0)}{\sum _{j > w_1}} a_{j1}t^j, t^{w_2} + \underset{j \notin \Gamma}{\sum _{j > w_2}} a_{j2}t^j , \cdots , t^{w_n} + \underset{j \notin \Gamma}{\sum _{j > w_n}} a_{jn}t^j \right).
\end{array}$$
\cqd

The elimination criteria presented in Lemma \ref{lemmaeliminasemigrupomaisv1menosv0esemigrupo} were originally considered for plane curves by Ebey and Zariski (see \cite{Ebey} and \cite{Zariski-book}). For a plane curve with semigroup $\Gamma=\langle v_0,v_1,\ldots ,v_g\rangle$, Zariski calls the corresponding form (\ref{curvagerallimpa}), namely  \begin{equation}
\label{Short-Param-Zariski}
\left( t^{v_0}, t^{v_1} + \underset{i \notin \Gamma \cup (\Gamma + v_1 - v_0)}{\sum _{i > v_1}} a_{1i}t^i \right),
\end{equation}a short parametrization.

\vspace{0.2cm}
In the last paragraph of Section \ref{Sec-Param-Curves}, we remark that any parametrization of a curve is $\mathcal{A}$-finitely determined. The previous result and \cite[Theorem 1.2]{Wall} give us an alternative proof of that result.

\begin{corollary}\label{corolariolimpezageral}
For any parametrization $\phi(t)=(t^{v_0},y_1(t),\ldots ,y_n(t))\in\mathcal{P}_{(v_0,w_1,\ldots ,w_n)}$ with semigroup $\Gamma$ and conductor $c$ we have $$
\phi(t) \underset{\mathcal{A}}{\sim} (t^{v_0},j^{e_1}y_1(t),\ldots ,j^{e_n}y_n(t))\ \ \mbox{where}\ \ e_i = \max \lbrace w_i, c-1 \rbrace.$$ In particular, $
\phi(t) \underset{\mathcal{A}}{\sim} j^{e}\phi(t)$, where $e=\max\{e_1,\ldots, e_n\}$.
\end{corollary}
\Dem It is sufficient to notice that $e_i+1=\max\{w_i+1,c\}=d_i$, where $d_i$ is defined in the proof of Theorem \ref{limpezageral}.\cqd

\begin{example}\label{example}
Let $\phi(t) \in\mathcal{P}_{(m,w_1,\ldots ,w_{m-1})}$ be a non-degenerate parametrized curve given by $$\phi(t)=\left(t^m,t^{w_1}+\sum_{j>w_1}a_{1j}t^j,\ldots ,t^{w_{m-1}}+\sum_{j>w_{m-1}}a_{m-1,j}t^j  \right).$$ 

Remark \ref{observacaodimensaoegenero} implies that the semigroup of $\phi$ is given by $\Gamma=\langle m,w_1,\ldots ,w_{m-1}\rangle$, while according to \cite{semigrupo}, its conductor is $c=w_{m-1}-m+1$. 

If $w_1\geq c$, then $\phi$ is $\mathcal{A}$-equivalent to 
$$\psi(t)=(t^m,t^{w_1},\ldots ,t^{w_{m-1}}).$$
In fact, since $m<w_1<\ldots <w_{m-1}$, if $w_1\geq c$ (or equivalently $c-1=w_{m-1}-w_1<m$), then $e_i=\max\{w_i,c-1\}=w_i$ for each $1\leq i<m$. The result is consequence of Corollary \ref{corolariolimpezageral}.

In particular, consider $$\phi (t) = \left(t^3, t^{7} + \sum_{j>7} a_{1j}t^j, t^{w_2} + \sum_{j>w_2} a_{2j}t^j \right) \in \mathcal{P}_{(3,7,w_2)}$$ with $w_2 \notin \langle 3, 7\rangle$ and $w_2 > 7$, that is, $w_2\in \{8,11\}$. In this case, we get $\Gamma=\langle 3, 7, w_2\rangle$ and its conductor is $c=w_2-2$. 

For $w_2=8$ we have that $\phi(t)$ is $\mathcal{A}$-equivalent to $\psi(t)=(t^3, t^7, t^8)$, since $c=6 < 7=w_1$.

On the other hand, if $w_2=11$, then we obtain $\Gamma=\langle 3,7,11\rangle$ with conductor $c=w_2-2=9$. Since $8\not\in\Gamma\cup (\Gamma+v_1-v_0)=\Gamma\cup (\Gamma+4)$, $e_1=\max\{w_1,c-1\}=8$ and $e_2=\max\{w_2,c-1\}=11$, by Corollary \ref{corolariolimpezageral}, we have that $\phi(t)$ is $\mathcal{A}$-equivalent to $\psi(t)=(t^3,t^7+at^8,t^{11})$ for some $a$.
\end{example}

Lemma \ref{lemmaeliminasemigrupomaisv1menosv0esemigrupo} establishes criteria to eliminate the terms $a_{ik}t^k$ from the component $y_i(t)$ of $\phi(t)\in\mathcal{P}_{(v_0,w_1,\ldots ,w_n)}$ whenever $k\in\Gamma$ or $k\in\Gamma+w_i-v_0$. In the next result, we show that if $y_i(t)\neq t^{w_i}$, we can eliminate all the terms in $y_i(t)$ with exponent $k+rv_0$, for any $k=v(y_i(t)- t^{w_i})$ and $r \geq 1$.

\begin{proposition}\label{propomega0}
Let $\phi(t)=(t^{v_0},y_1(t),\ldots ,y_n(t))\in\mathcal{P}_{(v_0,w_1,\ldots ,w_n)}$ be a parametrization as in \eqref{parametgeral2}, where $y_i(t)=t^{w_i}+\sum_{j>w_i}a_{ij}t^{j}$ for  $i=1,\ldots ,n$. If $y_i(t)-t^{w_i}=a_{ik}t^{k}+h.o.t.$ with $a_{ik}\neq 0$ for some $1\leq i\leq n$, then we can eliminate all the terms of order $k+rv_0$ in $y_i(t)$ with $r\geq 1$.
\end{proposition}
\Dem
For any $r\geq 1$ and some $\alpha \in \mathbb{C}$, let 
$$h_0=\alpha X^{r+1},\ \ h_i=\alpha\frac{w_i}{v_0}X^rY_i+Y_i\displaystyle \sum _{l = 2}^{\infty} \left( \begin{array}{c}
\frac{w_i}{v_0} \\
l
\end{array} \right) X^{lr},\ \mbox{and}\ h_s=0\ \ \mbox{for}\ 1\leq s\leq n,\ s\neq i.$$ Consider the pair $(\rho,\sigma) \in \mathcal{A}$ as in Proposition \ref{caracdiffeos}, and $\theta_s(t)$ as in (\ref{theta-i}) for $1 \leq s \leq n$. 

By Proposition \ref{acaodosdiffeo}, we obtain
$$\theta_i(t)=\frac{\alpha a_{ik}(w_i-k)}{v_0}t^{k+rv_0}+h.o.t.\hspace{0.3cm} \mbox{and} \hspace{0.3cm} \theta_s(t)=-\frac{\alpha w_s}{v_0}t^{w_s+rv_0} + h.o.t. ,$$ for every $1\leq s\leq n$ with $s\neq i$. Thus, $ord_t(\theta_s(t))=w_s+rv_0$, and $ord_t(\theta_i(t))=k+rv_0$. Taking $\alpha=-\frac{a_{i,k+rv_0}v_0}{a_{ik}(w_i-k)}$, by Remark \ref{remarkeliminacao}, we eliminate the term of order $k+rv_0$ in $y_i(t)$ and we obtain that $\phi(t)$ is $\mathcal{A}$-equivalent to $\bar{\phi}(t)=(t^{v_0},\bar{y}_1(t),\ldots ,\bar{y}_n(t))$ satisfying $$j^{k+rv_0}(\bar{y}_i(t))=j^{k+rv_0-1}(y_i(t)) \hspace{0.3cm} \mbox{and}\hspace{0.3cm} j^{w_s+rv_0-1}(\bar{y}_s(t))=j^{w_s+rv_0-1}(y_s(t))$$  for every $1\leq s\leq n$ with $s\neq i$.
\cqd

Zariski \cite[Lemma 5]{Zariski} shows that for a parametrized plane curve as in (\ref{Short-Param-Zariski}), if there exists $\lambda = \min \lbrace i\not\in\Gamma\cup (\Gamma+v_1-v_0) \,; \ a_{1i} \neq 0 \rbrace$, then $\lambda$ is invariant under $\mathcal{A}$-equivalence. In particular, if its semigroup $\Gamma=\langle v_0,v_1,v_2, \ldots, v_g \rangle$ satisfies $g\geq 2$, then $\lambda <v_2$.

For curves in $(\mathbb{C}^{n+1},0)$ with $n\geq 2$ admitting a parametrization $(t^{v_0},y_1(t),\ldots ,y_n(t))$ as (\ref{curvagerallimpa}), if there exists $\lambda = \min \lbrace i  \not\in\Gamma\cup (\Gamma+v_1-v_0); \ a_{1i} \neq 0 \rbrace$, by Proposition \ref{propomega0}, we can eliminate every term in $y_1(t)$ of order $\lambda+rv_0$ for $r\geq 1$. 
Although we cannot conclude that $\lambda <v_2$ in this higher-dimensional setting, under this assumption $\lambda$ is also an invariant with respect to the action of diffeomorphisms given in Proposition \ref{caracdiffeos}.

\begin{theorem}\label{lambdazariski}
Let $$\phi (t) = \displaystyle  \left (t^{v_0}, t^{w_1}+\underset{j \notin \Gamma \cup (\Gamma + v_1 - v_0)}{\sum _{j > v_1}} a_{1j}t^j, y_2(t), \ldots, y_n(t)\right )\in\mathcal{P}_{(v_0,w_1,\ldots ,w_n)}$$ be a parametrization as \eqref{curvagerallimpa}, with semigroup $\Gamma = \langle v_0, v_1, v_2, \ldots, v_g \rangle$ and suppose that there exists $\lambda = min \lbrace j\not\in\Gamma\cup (\Gamma+v_1-v_0) \,; \ a_{1j} \neq 0 \rbrace$. 
If $\lambda < v_2$, then $\lambda$ is invariant with respect to the action of diffeomorphisms given in Proposition \ref{caracdiffeos}.
\end{theorem}
\Dem Let $\psi (t) = \displaystyle \left( t^{v_0}, t^{v_1} + b_{1\gamma}t^{\gamma} + \sum _{j > \gamma} b_{1j}t^j, \overline{y}_2(t), \ldots, \overline{y}_n(t) \right) \in \mathcal{P}_{(v_0, w_1, \ldots, w_n)}$ as in  \eqref{curvagerallimpa}, with $v_1=w_1$, and $\gamma=\min\{j\not\in\Gamma\cup (\Gamma+v_1-v_0);\ b_{1j}\neq 0\}$. In Proposition \ref{caracdiffeos}, we proved that if $(\rho, \sigma) \in \mathcal{A}$ is such that $(\sigma \circ \phi \circ \rho ^{-1})(t) = \psi (t) \in \mathcal{P}_{(v_0, w_1, \ldots, w_n)}$, then
$\rho (t) = t \cdot \left( \epsilon ^{v_0} + \frac{h _0(\phi (t))}{t^{v_0}} \right)^{\frac{1}{v_0}}$ and $\sigma = (\sigma _0, \sigma _1, \ldots, \sigma _n)$, with $$
\sigma _0 (X, \underline{Y})  =  \epsilon ^{v_0}X + h _0(X, \underline{Y}) \ \ \ \mbox{and}\ \ \
\sigma _k(X, \underline{Y})  =  \epsilon ^{w_k} Y_k + h _k (X, \underline{Y}),
$$
where $\epsilon \in \mathbb{C} ^*$ and $h _k \in \langle Y_{k + 1}, \ldots, Y_n \rangle + \mathcal{M}_{n + 1}^2$ for $0\leq k\leq n$ such that $ord _t(h _k(\phi (t))) > w_k$, for $k \in \lbrace 1, \ldots, n \rbrace$. Without loss of generality, we may assume that $\epsilon = 1$.

The condition $(\sigma \circ \phi \circ \rho ^{-1})(t) = \psi (t)$ implies that $(\sigma \circ \phi)(t) = (\psi \circ \rho)(t)$. Consequently, its second component satisfies \begin{equation}\label{igualdadelambdazariski}
y_1(t)+h_1(\phi(t))=t^{v_1} + a_{1\lambda }t^{\lambda} + h _1(\phi (t)) +h.o.t.= (\rho (t))^{v_1} + b_{1\gamma }(\rho (t))^{\gamma} + h.o.t.
\end{equation}

Suppose that $\lambda < \gamma$.

Notice that there is no term of order $\lambda$ in the expansion of $h_1(\phi(t))$. Indeed, since $h_1(X,\underline{Y})\in\langle Y_2,\ldots ,Y_n\rangle+\mathcal{M}^2_{n+1}$ and $\lambda<v_2\leq w_2=\min (\Gamma\setminus\Gamma_2)<w_3<\ldots <w_n$ (see (\ref{gammai})), the only possibility of obtaining a term of order $\lambda$ in $h_1(\phi(t))$ should be for $h_1 \in \langle X,Y_1 \rangle ^2$. It is sufficient to assume that $h_1(X,Y_1)$ is a monomial in $\langle X, Y_1\rangle^2$. However, if $h_1(X,Y_1)=X^aY_1^b$ with $a+b\geq 2$ we have
$$h_1(\phi(t))=t^{av_0}(t^{v_1}+a_{1\lambda}t^{\lambda}+h.o.t.)^b=t^{av_0+bv_1}+ba_{1\lambda}t^{\lambda+(b-1)v_1+av_0}+h.o.t.,$$
with $av_0+bv_1\in\Gamma$ and $\lambda+(b-1)v_1+av_0>\lambda$. So, by (\ref{igualdadelambdazariski}), the only term of order $\lambda$ in $(\sigma\circ\phi)(t)$ is $a_{1\lambda}t^{\lambda}$ and every term with order $\gamma<\lambda$ is such that $\gamma\in\Gamma$.

On the other hand, since $\lambda < \gamma$, the term $a_{1\lambda }t^{\lambda}$ must appear in the expansion of
$$
(\rho (t))^{v_1} = t^{v_1} \left( 1 + \frac{h _0(\phi (t))}{t^{v_0}} \right)^{\frac{v_1}{v_0}} = t^{v_1} +  \sum _{l = 1}^{\infty} \left(
\begin{array}{c}
\frac{v_1}{v_0} \\
l
\end{array}
\right) t^{v_1-lv_0}h _0^l(\phi (t)).
$$
More precisely, $\lambda$ corresponds to the term of lowest order that does not belong to $\Gamma\cup (\Gamma+v_1-v_0)$ in $t^{v_1-lv_0}h_0^l(\phi(t))$ for some $l\geq 1$.
Since $\lambda < v_2\leq w_2$, it is sufficient to consider $h _0$ as a monomial in $\langle Y_1 \rangle + \langle X, Y_1 \rangle ^2$, that is, $X^aY_1^b$ or $Y_1^d$ with $a+b \geq 2$ and $d\geq 1$.

If $h _0(X, \underline{Y}) = X^aY_1^b$ with $a,b\geq 1$ then, for any $l\geq 1$, we get
$$
t^{v_1-lv_0}h_0^l(\phi(t))=t^{v_1+lv_0(a-1)}\cdot \left (t^{lbv_1}+lba_{1\lambda}t^{\lambda+(lb-1)v_1}+h.o.t.\right ).$$
Since $v_1+lv_0(a-1)+lbv_1\in\Gamma$ and $v_1+lv_0(a-1)+\lambda+(lb-1)v_1>\lambda$, in this case, no term with the order $\lambda$ appears.

If $h_0(X,\underline{Y})=Y_1^d$ for some $d\geq 1$, then
$$
t^{v_1-lv_0}h_0^l(\phi(t))=t^{v_1-lv_0}\cdot \left (t^{ldv_1}+lda_{1\lambda}t^{\lambda+(ld-1)v_1}+h.o.t.\right ).$$
Note that $v_1-lv_0+\lambda+(ld-1)v_1=\lambda+l(dv_1-v_0)>\lambda$ because $l,d\geq 1$ and $v_1>v_0$. So, the only possibility to obtain a term of order $\lambda$ in $(\psi\circ\rho)(t)$ is $\lambda=(ld+1)v_1-lv_0$ for some $d,l\geq 1$.

Since $d,l\geq 1$ and $v_1>v_0$, we cannot have $\lambda <(d+1)v_1- v_0$. In addition, if $\lambda=(d+1)v_1-v_0=dv_1+v_1-v_0$, then $\lambda\in\Gamma+v_1-v_0$, which is not the case. Hence, we must have $$
(d+1)v_1 - v_0< \lambda
\  \ \ \mbox{for some}\ \ \ d \geq 1.$$
In this case, $(d+1)v_1-v_0\in \Gamma$, otherwise, the right side of \eqref{igualdadelambdazariski} would contain a term of order less than $\lambda$ that does not belong to $\Gamma$, which does not occur on the left side of that equation. 

As $(d+1)v_1-v_0\in \Gamma$
 and $(d+1)v_1-v_0<\lambda<v_2$, it follows that $(d + 1)v_1 - v_0 \in\langle v_0,v_1\rangle$, that is, there exist $\alpha,\beta\in\mathbb{N}$ such that $(d+1)v_1-v_0= \alpha v_1 + \beta v_0$. Thus, \begin{equation}
\label{Relation-lambda}
(d+ 1 - \alpha)v_1 = (\beta + 1)v_0.
\end{equation} However, $v_1\not\in\langle v_0\rangle$ and $v_0<v_1$; so, we must have $\alpha<b$. Thus, $d+1 - \alpha >1$ and $\beta>1$.

Now, we write $l = 2q + r$, where $q, r \in \mathbb{N}$ with $r \in \lbrace 0, 1 \rbrace$. Therefore, \begin{equation}\label{auxlambdazariski}
\lambda=l(d v_1 - v_0) + v_1=  q(2dv_1 - 2v_0) + (dr + 1)v_1 - rv_0
\end{equation}
and $$
(dr + 1)v_1 - rv_0 = \left\lbrace 
\begin{array}{cc}
    dv_1 + v_1 - v_0, & \mbox{if} \ r = 1 \\
    v_1, & \mbox{if} \ r = 0 ,
\end{array}
\right. 
$$
that is, $(dr + 1)v_1 - rv_0 \in \Gamma \cup (\Gamma + v_1 - v_0)$.

By (\ref{Relation-lambda}), we obtain the equality $(d + 1)v_1 - v_0 = \alpha v_1 + \beta v_0$, with $d-1+\alpha,\beta-1>1$. Thus, $2dv_1 - 2v_0= (d+1)v_1-v_0+dv_1-v_1-v_0=(d - 1 + \alpha)v_1 + (\beta - 1)v_0 \in \Gamma$. Therefore, by \eqref{auxlambdazariski}, we conclude that $$
\lambda=l(dv_1 - v_0) + v_1 = q(2dv_1 - 2v_0) + (dr + 1)v_1 - rv_0 \in \Gamma \cup (\Gamma + v_1 - v_0),
$$
which is a contradiction, since $\lambda\not\in\Gamma\cup (\Gamma+v_1-v_0)$. Consequently, in all the analyzed cases, if $(\sigma\circ\phi \circ \rho^{-1})(t)=\psi(t)$, then we cannot have $\lambda<\gamma$. Applying the same argument to $(\phi \circ \rho ^{-1})(t) = (\sigma ^{-1} \circ \psi)(t)$, we conclude that $\gamma < \lambda$ does not occur. 

Hence, if $(\sigma\circ\phi \circ \rho ^{-1})(t) = \psi(t)$, then we must have $\lambda = \gamma$.
\cqd

Notice that, in Example \ref{example}, we found that any parametrization with values semigroup $\Gamma=\langle 3,7,11\rangle$ is $\mathcal{A}$-equivalent to $(t^3,t^7+at^8,t^{11})$. Since $8\not\in\Gamma\cup (\Gamma+v_1-v_0)$ and $8<v_2=11$, if $a\neq 0$, then Theorem \ref{lambdazariski} gives us that $8$ is an $\mathcal{A}$-invariant.

The hypothesis $\lambda < v_2$ in Theorem \ref{lambdazariski} cannot be removed as we illustrate in the following example: 

\begin{example}\label{exemplo-lambda}
Let us consider $\phi (t) = (t^9, t^{12} + t^{32}, t^{19}, t^{22})\in\mathcal{P}_{(9,12,19,22)}$, whose values semigroup is $$
\Gamma = \lbrace 0, 9, 12, 18, 19, 21, 22, 24, 27, 28, 30, 31, 33, 34, 36, \ldots \rbrace = \langle 9, 12, 19, 22 \rangle 
$$ with conductor $c = 36$. 
Note that $\lambda = 32$ satisfies $\lambda > 19 = w_2 = v_2$ and $32 \notin \Gamma \cup (\Gamma + v_1 - v_0)$.
Taking $\epsilon=1$, $h _0 = \alpha Y_2$, $h _1 =  \frac{12}{9}\alpha Y_3$, $h _2 = h _3 = 0$, for some $\alpha \in \mathbb{C}$, and considering $(\rho,\sigma)\in\tilde{\mathcal{A}}_1$ as given in Proposition \ref{caracdiffeos}, we obtain by Proposition \ref{acaodosdiffeo} that 
$$
\begin{array}{rcl}
\theta _1(t) & = & \displaystyle - \frac{32}{9}\alpha t^{42} - t^{12} \sum _{l = 2}^{\infty} \left(
\begin{array}{c}
\frac{12}{9} \\
l
\end{array}
\right) \left( \alpha t^{10} \right)^l -  t^{32}\sum _{l = 2}^{\infty} \left( 
\begin{array}{c}
\frac{32}{9} \\
l
\end{array}
\right) \left( \alpha t^{10} \right)^l.
\end{array}
$$
Thus, $ord _t(\theta _1(t)) = \displaystyle ord _t \left(  t^{12} \sum _{l = 2}^{\infty} \left(
\begin{array}{c}
\frac{12}{9} \\
l
\end{array}
\right) \left( \alpha t^{10} \right)^l \right) = 32<ord _t(\theta _i(t))=19+w_i-9$ for $i \in \lbrace 2, 3 \rbrace$.
By Remark \ref{remarkeliminacao}, we can eliminate the term $t^{32}$ from the component $y_1(t)=t^{12}+t^{32}$ by choosing an appropriate value for $\alpha\in\mathbb{C}$. 
\end{example}

Until now, we have considered the $\tilde{\mathcal{A}}_1$-action on a parametrized curve $\phi$ to eliminate some terms, since the group $\mathcal{H}$ neither introduces nor eliminates any term in $\phi$ (see \eqref{grupoH}). However, we can use such an action to normalize a nonzero coefficient in a parametrization as we show in the next proposition.

\begin{proposition}\label{normalizaocoeficiente}
Given $\phi (t) = \left( t^{v_0},t^{w_1} + \sum _{j > w_1}a_{1j}t^j,  \cdots, t^{w_n} + \sum _{j > w_n}a_{nj}t^j \right)$ with $a_{sk} \neq 0$ for some $1\leq s\leq n$ and $k>w_s$, there exists $(\rho, \sigma) \in \mathcal{H}$ such that $$ (\sigma \circ \phi \circ \rho ^{-1})(t) = \left( t^{v_0}, t^{w_1} +  \sum _{j > w_1}b_{1j}t^j, \cdots, t^{w_n} + \sum _{j > w_n}b_{nj}t^j \right)$$ where $b_{ij}=a_{ij}\cdot a_{sk}^{ \frac{w_i-j}{k - w_s}}$. In particular, $b_{sk}=1$.
\end{proposition}
\Dem
Considering $
\rho ^{-1}(t) = a_{sk}^{- \frac{1}{k - w_s}}
t$ and $\sigma(X,\underline{Y})=(\sigma_0(X,\underline{Y}),\ldots ,\sigma_n(X,\underline{Y}))$ with 
$$\sigma_0(X,\underline{Y})=a_{sk}^{ \frac{v_0}{k - w_s}}
X\ \ \ \mbox{and}\ \ \ \sigma_i(X,\underline{Y})=a_{sk}^{ \frac{w_i}{k - w_s}}Y_i\ \ \ \mbox{for every}\ \ 1\leq i\leq n,
$$
we get 
$(\sigma\circ\phi\circ\rho^{-1})(t)=(\sigma_0\circ\phi\circ\rho^{-1}(t),\ldots ,\sigma_n\circ\phi\circ\rho^{-1}(t))$, where
$$(\sigma_0\circ\phi\circ\rho^{-1})(t)=a_{sk}^{\frac{v_0}{k - w_s}}\cdot a_{sk}^{- \frac{v_0}{k-w_s}}t^{v_0}=t^{v_0}\ \ \ \mbox{and}$$
$$(\sigma_i\circ\phi\circ\rho^{-1})(t)=a_{sk}^{ \frac{w_i}{k - w_s}}\left ( a_{sk}^{- \frac{w_i}{k - w_s}}
t^{w_i}+\sum_{j>w_i}a_{ij}\cdot a_{sk}^{- \frac{j}{k - w_s}}
t^j\right )=t^{w_i}+\sum_{j>w_i}a_{ij}\cdot a_{sk}^{ \frac{w_i-j}{k - w_s}}
t^j.$$ Then, $(\sigma_i\circ\phi\circ\rho^{-1})(t)=t^{w_i}+\sum_{j>w_i}b_{ij}
t^j,$ where $b_{ij}=a_{ij}\cdot a_{sk}^{ \frac{w_i-j}{k - w_s}}$. In particular, $b_{sk}=1$.
\cqd

\begin{example}\label{ex-3}
    Let $\phi(t)=(t^3,t^{w_1}+\sum_{j>w_1}a_{1j},t^{w_2}+\sum_{j>w_2}a_{2j}t^j)$ be a non-degenerate parametrized curve as (\ref{parametgeral2}). 
    
    By Remark \ref{observacaodimensaoegenero}, it follows that $\Gamma=\langle 3,v_1,v_2\rangle$, where $v_1=w_1$ and $v_2=w_2$. In Section \ref{Sec-Semigroups}, we note that the Ap\'ery set of $\Gamma$ is $A=\{0,v_1,v_2\}$, its conductor is $v_2-2$ and the gaps set of $\Gamma$ is $$\mathbb{N}\setminus\Gamma=\left \{v_1-3i;\ 0<i<\frac{v_1}{3}\right \}\cup\left \{v_2-3i;\ 0<i<\frac{v_2}{3}\right \}.$$
Since $\{\gamma\in\mathbb{N};\ \gamma>v_1\ \mbox{and}\ \gamma\not\in\Gamma\cup (\Gamma+v_1-3)\}=\{v_2-3j;\ 1\leq j<\frac{v_2-v_1}{3}\}$, by Theorem \ref{limpezageral}, we have that $\phi(t)$ is $\mathcal{A}$-equivalent to
$$\phi_1(t)=\left (t^3,t^{v_1}+\sum_{1\leq j<\frac{v_2-v_1}{3}}b_{j}t^{v_2-3j},t^{v_2}\right ).$$

If $b_{j}=0$ for every $1\leq j<\frac{v_2-v_1}{3}$, then $\phi_1(t)=(t^3,t^{v_1},t^{v_2})$.

If there exists $b_j\neq 0$ for some $1\leq j<\frac{v_2-v_1}{3}$, then we consider $$k:=\max \left\{j;\ b_{j}\neq 0\ \mbox{and}\ 1\leq j<\frac{v_2-v_1}{3}\right\}.$$ As $v_2-3k<v_2$, by Theorem  \ref{lambdazariski}, we have that $v_2-3k$ is an $\tilde{\mathcal{A}}_1$-invariant. In addition, by Proposition \ref{propomega0}, we can eliminate in the second component of $\phi_1(t)$ all terms of the form $v_2-3k+3r$, for every $r\geq 1$, that is, all terms greater than $v_2-3k$. Therefore, by Proposition \ref{normalizaocoeficiente}, $\phi_1(t)$ (consequently $\phi(t)$) is $\mathcal{A}$-equivalent to $(t^3,t^{v_1}+t^{v_2-3k},t^{v_2})$.

Hence, any non-degenerate parametrized curve of multiplicity $3$ in $\mathbb{C}^3$  is $\mathcal{A}$-equivalent to \begin{equation}\label{normal-form-3}(t^3,t^{v_1},t^{v_2})\ \ \mbox{or}\ \ (t^3,t^{v_1}+t^{v_2-3k},t^{v_2}),
\end{equation}
where $3<v_1<v_2$, $gcd(3,v_1,v_2)=1$ and $1\leq k<\frac{v_2-v_1}{3}$. This result was also obtained by Gibson and Hobbs in \cite[Proposition 6.3]{GibsonHobbs}.
\end{example}

Any plane branch with multiplicity $3$ admits the numerical semigroup $\Gamma=\langle 3,v_1\rangle$ for some $v_1>3$ with $\gcd(3,v_1)=1$. In this case, the Ap\'ery set of $\Gamma$ is $\{0,v_1,2v_1\}$ and its conductor is $2v_1-2$. 
According to Zariski \cite{Zariski-book}, any such plane branch is analytically equivalent to one with the following parametrizations: $$(t^3,t^{v_1}) \hspace{0.3cm} \mbox{or} \hspace{0.3cm}(t^3,t^{v_1}+t^{2v_1-3j}),$$ for some integer $j$ satisfying $2\leq j<\frac{v_1}{3}$. Using the previous plane normal forms, we can obtain the normal forms for branches with multiplicity $3$ in $\mathbb{C}^3$ presented in (\ref{normal-form-3}). In fact, let $v_2\not\in\langle 3,v_1\rangle$ be an integer such that $v_2>v_1$ and $v_2\equiv 2v_1\mod 3$. This allows us to write $v_2=2v_1-3i$ for some $1\leq i< \frac{v_1}{3}$. Now, consider the parametrizations in $\mathbb{C}^3$ of the form: $$(t^3,t^{v_1},t^{v_2}) \hspace{0.3cm} \mbox{or} \hspace{0.3cm}(t^3,t^{v_1}+t^{2v_1-3j},t^{v_2}),$$ where $2v_1-3j\leq v_2-3,$ for $2\leq j<\frac{v_1}{3}$. Since $v_2=2v_1-3i$ and $2v_1-3j\leq v_2-3$ we get $j\geq i+1$. Hence, setting $k=j-i$, we obtain $2v_1-3j=v_2-3k$ with $1\leq k< \frac{v_2-v_1}{3}$. Thus, we recover the parametrizations in (\ref{normal-form-3}) that naturally arise from the plane ones. 

The construction illustrated above will be generalized in the next theorem.

Let $\phi_0 (t) = (t^{v_0}, y_1(t), \ldots, y_n(t))\in\mathcal{P}_{(v_0,w_1,\ldots ,w_n)}$ be a non-degenerate parametrization with semigroup $\Gamma_0 = \langle v_0, \ldots, v_g \rangle$ and Ap\'ery set $\{0,a_1,\ldots ,a_{v_0-1}\}$. Given $v_{g+1}\in\mathbb{N}\setminus\Gamma_0$ satisfying $v_{g+1}\equiv a_{v_0-1}\mod v_0$ and $\max\{v_g,a_{v_0-2}\}<v_{g+1}<a_{v_0-1}$ consider a parametrization
$\phi_1(t) = (t^{v_0},y_1(t),\ldots ,y_n(t),t^{v_{g+1}}+\sum_{j>v_{g+1}}a_jt^j)\in\mathcal{P}_{(v_0,w_1,\ldots ,w_n,v_{g+1})}$. It follows that $\phi_1$ admits the semigroup $\Gamma_1=\langle v_0,\ldots ,v_g,v_{g+1}\rangle$. By Proposition \ref{diminuircond}, the Ap\'ery set of $\Gamma_1$ is $\lbrace 0, a_1, \ldots, a_{v_0 - 2},v_{g+1} \rbrace$ and its conductor is $v_{g + 1} - v_0 + 1$.

\begin{theorem}\label{eliminadeumaiffdaoutra}
Let $\phi_0 (t) = (t^{v_0}, y_1(t), \ldots, y_n(t))\in\mathcal{P}_{(v_0,w_1,\ldots ,w_n)}$ be a non-degenerate parametrized curve and set $$\phi_1(t) = \left(t^{v_0},y_1(t),\ldots ,y_n(t),t^{v_{g+1}}+\sum_{j>v_{g+1}}a_jt^j\right)\in\mathcal{P}_{(v_0,w_1,\ldots ,w_n,v_{g+1})}$$ as described above. A term of order $l<v_{g + 1}$ can be eliminated from the component $y_i(t)$ in $\phi_0(t)$, through the action of an element in $\tilde{\mathcal{A}}_1$ as given in Proposition \ref{caracdiffeos}, if and only if it can be eliminated from $y_i(t)$ in $\phi_1(t)$, for $i \in \lbrace 1, \ldots, n \rbrace$. 
Moreover, $\phi_1(t)$ is $\mathcal{A}$-equivalent to $$(t^{v_0}, j^{k_1}y_1(t), j^{k_2}y_2(t), \ldots, j^{k_n}y_n(t), t^{v_{g + 1}})\in \mathcal{P}_{(v_0,w_1,\ldots ,w_n,v_{g+1})},$$
where $k_i = \max \lbrace w_i, v_{g + 1} - v_0 \rbrace$ for $1 \leq i \leq n$.
\end{theorem}
\Dem Suppose that for some $i \in \{1,\ldots, n\}$, we can eliminate a term of order $l<v_{g+1}$ from the component $y_i(t)$ of $\phi_0(t)\in\mathcal{P}_{(v_0,w_1,\ldots ,w_n)}$ by the action of an element $(\rho,\sigma)\in \tilde{\mathcal{A}}_1$ as characterized in Proposition \ref{caracdiffeos}, where $\sigma(X,\underline{Y})=(\sigma_0(X,\underline{Y}),\ldots ,\sigma_n(X,\underline{Y}))$. 

By Proposition \ref{caracdiffeos}, there exist $h_j \in  \langle Y_{j+1}, \ldots, Y_n \rangle + \mathcal{M}_{n + 1}^2$ for $j\in\{0,\ldots ,n\}$ such that $\rho(t)=t\cdot\left ( 1+\frac{h_0(\phi_0(t))}{t^{v_0}}\right )^{\frac{1}{v_0}}$, $\sigma_0(X,\underline{Y})=X+h_0(X,\underline{Y})$ and $\sigma_k(X,\underline{Y})=Y_k+h_k(X,\underline{Y})$ for $1\leq k\leq n$. Since $\langle Y_{j + 1}, \ldots, Y_n \rangle + \mathcal{M}_ {n + 1}^2 \subset \langle Y_{j + 1}, \ldots, Y_n, Y_{n + 1} \rangle +  \mathcal{M}_{n + 2}^2$, taking $(\rho,\sigma')$ with $$\sigma'(X,\underline{Y},Y_{n+1})=(\sigma_0(X,\underline{Y}),\ldots ,\sigma_n(X,\underline{Y}),Y_{n+1}),$$ we eliminate the same term of order $l$ from the component $y_i(t)$ of $\phi_1(t)$.

Conversely, suppose that for some $i\in \{1,\ldots, n\}$, a term of order $l<v_{g + 1}$ can be eliminated from the component $y_i(t)$ of $\phi_1(t)$ via an element $(\rho,\sigma)\in \tilde{\mathcal{A}}_1$ given by $$\rho(t)=t\cdot \left ( 1+\frac{h_0(\phi_1(t))}{t^{v_0}}\right )\ \ \ \mbox{and}\ \ \ \sigma(X,\underline{Y},Y_{n+1})=(\sigma_0(X,\underline{Y},Y_{n+1}),\ldots , \sigma_{n+1}(X,\underline{Y},Y_{n+1}))$$ with $\underline{Y}=(Y_1,\ldots,Y_n)$, $\sigma_0(X,\underline{Y},Y_{n+1})=X+h_0(X,\underline{Y},Y_{n+1})$, $\sigma_k(X,\underline{Y},Y_{n+1})=Y_k+h_k(X,\underline{Y},Y_{n+1})$ for $1\leq k\leq n+1$ and $h_j\in\langle Y_{j+1},\ldots ,Y_{n+1}\rangle +\mathcal{M}^2_{n+2}$ for $0\leq j\leq n+1$.

Writing $y_{n+1}(t)=t^{v_{g+1}}+\sum_{j>v_{g+1}}a_jt^j$, the condition $l<v_{g+1}=ord_t(y_{n+1}(t))$ implies that $h_j \in  \langle Y_{j+1}, \ldots, Y_n \rangle + \mathcal{M}_{n + 1}^2$, and $\sigma_j(X,\underline{Y},Y_{n+1})=\sigma''_j(X,\underline{Y})$ for $0\leq j\leq n$. Consequently, we have $h_j(\phi_1(t))=h_j(\phi_0(t))$, which means that the element $(\rho,\sigma'')$ defined by $\sigma''(X,\underline{Y})=(\sigma''_0(X,\underline{Y}),\ldots ,\sigma''_{n}(X,\underline{Y}))$ eliminates the same term of order $l$ from $y_i(t)$ of $\phi_0(t)$.

Moreover, as $\phi_1(t)$ admits the semigroup $\Gamma_1= \langle v_0, v_1, \ldots, v_g, v_{g + 1} \rangle$, with conductor $v_{g + 1} -v_0+1$, any integer $l>v_{g+1}-v_0$ belongs to $\Gamma_1$. Therefore, if there exists a term of order $l>v_{g+1}-v_0$ in $y_i(t)$, it can be eliminated as indicated in item 1 of Lemma \ref{lemmaeliminasemigrupomaisv1menosv0esemigrupo}. The same argument to the component $y_{n+1}(t)=t^{v_{g+1}}+\sum_{j>v_{g+1}}a_jt^j$ allows us to conclude that $\phi_1(t)$ is $\mathcal{A}$-equivalent to $$\left( t^{v_0}, j^{k_1}(y_1(t)), \ldots, j^{k_n}(y_n(t)), t^{v_{g + 1}} \right)\in\mathcal{P}_{(v_0,w_1,\ldots ,w_n,v_{g+1})},$$
where $k_i = \max \lbrace w_i, v_{g + 1} - v_0 \rbrace$, for $1 \leq i \leq n$.
\cqd 

The results of this paper will be applied to construct normal forms for parametrized curves of multiplicity 4 in $\mathbb{C}^3$ and $\mathbb{C}^4$ (as presented in \cite{HRoR}).

\end{document}